\documentclass{icmart}

\usepackage{amsmath}
\usepackage{amssymb}
\usepackage{amscd}
\usepackage{mathrsfs}
\usepackage{pgfkeys}
\usepackage{pgfopts}
\usepackage{xcolor}
\usepackage{ytableau}

\usepackage{xy}
\xyoption{all}

\newcommand{\nc}{\newcommand}

\nc{\kk}{{\mathsf{k}}}

\nc{\SB}{{\mathsf{B}}}
\nc{\HH}{{\mathsf{HH}}}
\nc{\SJ}{{\mathsf{J}}}
\nc{\SK}{{\mathsf{K}}}
\nc{\SL}{{\mathsf{L}}}
\nc{\SM}{{\mathsf{M}}}
\nc{\SO}{{\mathsf{O}}}
\nc{\SQ}{{\mathsf{Q}}}
\nc{\SR}{{\mathsf{R}}}
\nc{\SSS}{{\mathsf{S}}}
\nc{\ST}{{\mathsf{T	}}}
\nc{\SU}{{\mathsf{U}}}

\nc{\C}{{\mathbb{C}}}
\nc{\LL}{{\mathbb{L}}}
\nc{\PP}{{\mathbb{P}}}
\nc{\QQ}{{\mathbb{Q}}}
\nc{\RR}{{\mathbb{R}}}
\nc{\SSSS}{{\mathbb{S}}}
\nc{\ZZ}{{\mathbb{Z}}}

\nc{\ba}{{\mathbf{a}}}
\nc{\bm}{{\mathbf{m}}}
\nc{\bq}{{\mathbf{q}}}
\nc{\bu}{{\mathbf{u}}}
\nc{\BB}{{\mathbf{B}}}
\nc{\BC}{{\mathbf{C}}}
\nc{\BD}{{\mathbf{D}}}
\nc{\BG}{{\mathbf{G}}}
\nc{\BH}{{\mathbf{H}}}
\nc{\BK}{{\mathbf{K}}}
\nc{\BL}{{\mathbf{L}}}
\nc{\BM}{{\mathbf{M}}}
\nc{\BP}{{\mathbf{P}}}
\nc{\BS}{{\mathbf{S}}}
\nc{\BT}{{\mathbf{T}}}
\nc{\BU}{{\mathbf{U}}}
\nc{\BZ}{{\mathbf{Z}}}
\nc{\BPr}{{\mathsf{P}}}
\nc{\BR}{{\mathbf{R}}}
\nc{\BW}{{\mathbf{W}}}

\nc{\CA}{{\mathscr{A}}}
\nc{\CB}{{\mathscr{B}}}
\nc{\TCB}{\tilde{\mathscr{B}}}
\nc{\CC}{{\mathscr{C}}}
\nc{\TCC}{\tilde{\mathscr{C}}}
\nc{\D}{{\mathcal{D}}}
\nc{\CE}{{\mathcal{E}}}
\nc{\CF}{{\mathcal{F}}}
\nc{\CG}{{\mathcal{G}}}
\nc{\CH}{{\mathcal{H}}}
\nc{\CI}{{\mathcal{I}}}
\nc{\CJ}{{\mathscr{J}}}
\nc{\CK}{{\mathscr{K}}}
\nc{\CL}{{\mathscr{L}}}
\nc{\CM}{{\mathcal{M}}}
\nc{\CN}{{\mathcal{N}}}
\nc{\CO}{{\mathcal{O}}}
\nc{\CP}{{\mathcal{P}}}
\nc{\CQ}{{\mathcal{Q}}}
\nc{\CR}{{\mathcal{R}}}
\nc{\CS}{{\mathcal{S}}}
\nc{\CT}{{\mathcal{T}}}
\nc{\CU}{{\mathcal{U}}}
\nc{\CV}{{\mathcal{V}}}
\nc{\CW}{{\mathcal{W}}}
\nc{\CX}{{\mathscr{X}}}
\nc{\CY}{{\mathcal{Y}}}

\nc{\fa}{{\mathfrak{a}}}
\nc{\fb}{{\mathfrak{b}}}
\nc{\fd}{{\mathfrak{d}}}
\nc{\fg}{{\mathfrak{g}}}
\nc{\fn}{{\mathfrak{n}}}
\nc{\fp}{{\mathfrak{p}}}
\nc{\fu}{{\mathfrak{u}}}

\nc{\FA}{{\mathfrak{A}}}
\nc{\FB}{{\mathfrak{B}}}
\nc{\FD}{{\mathfrak{D}}}
\nc{\FE}{{\mathfrak{E}}}
\nc{\FL}{{\mathfrak{L}}}
\nc{\FM}{{\mathfrak{M}}}

\nc{\TY}{{\tilde{Y}}}

\nc{\lotimes}{\mathbin{\mathop{\otimes}\limits^{\mathbb{L}}}}

\nc{\CExt}{\mathop{\mathcal{E}\mathit{xt}}\nolimits}
\nc{\CHom}{\mathop{\mathcal{H}\mathit{om}}\nolimits}
\nc{\CEnd}{\mathop{\mathcal{E}\mathit{nd}}\nolimits}
\nc{\RCHom}{\mathop{\mathsf{R}\mathcal{H}\mathit{om}}\nolimits}

\nc{\Hom}{\mathop{\mathsf{Hom}}\nolimits}
\nc{\Ext}{\mathop{\mathsf{Ext}}\nolimits}
\nc{\End}{\mathop{\mathsf{End}}\nolimits}
\nc{\RGamma}{\mathop{{\mathsf{R}}\Gamma}\nolimits}
\nc{\RHom}{\mathop{\mathsf{RHom}}\nolimits}

\nc{\Tor}{\mathop{\mathsf{Tor}}\nolimits}

\nc{\Hilb}{\mathop{\mathsf{Hilb}}\nolimits}
\nc{\Spec}{\mathop{\mathsf{Spec}}\nolimits}
\nc{\Proj}{\mathop{\mathsf{Proj}}\nolimits}
\nc{\Pic}{\mathop{\mathsf{Pic}}\nolimits}
\nc{\Br}{\mathop{\mathsf{Br}}\nolimits}

\nc{\Mod}{\mathop{\mathsf{Mod}}\nolimits}
\nc{\grmod}{\mathop{\mathsf{grmod}}\nolimits}
\nc{\Grmod}{\mathop{\mathsf{Grmod}}\nolimits}
\nc{\Qcoh}{\mathop{\mathsf{Qcoh}}\nolimits}

\nc{\Ann}{\mathop{\mathsf{Ann}}\nolimits}
\nc{\Ker}{\mathop{\mathsf{Ker}}\nolimits}
\nc{\Coker}{\mathop{\mathsf{Coker}}\nolimits}

\nc{\Cone}{\mathop{\mathsf{Cone}}\nolimits}
\nc{\Tot}{{\mathsf{Tot}}}

\nc{\coh}{{\mathop{{\mathsf{coh}}}}}
\nc{\Ab}{{\mathop{\mathcal{A}\mathit{b}}}}

\nc{\Tr}{\mathop{\mathsf{Tr}}\nolimits}
\nc{\Ind}{\mathop{\mathsf{Ind}}\nolimits}
\nc{\Res}{\mathop{\mathsf{Res}}\nolimits}

\nc{\Conv}{\mathop{\mathsf{Conv}}\nolimits}

\nc{\codim}{\mathop{\mathsf{codim}}\nolimits}

\nc{\sing}{{\mathsf{sing}}}
\nc{\supp}{\mathop{\mathsf{supp}}}
\nc{\vol}{\mathop{\mathsf{vol}}\nolimits}
\nc{\ch}{\mathop{\mathsf{ch}}\nolimits}
\nc{\perf}{{\mathsf{perf}}}
\nc{\chr}{\mathop{\mathsf{char}}}
\nc{\rk}{\mathop{\mathsf{rk}}}
\nc{\Pf}{{\mathsf{Pf}}}
\nc{\Gr}{{\mathsf{Gr}}}
\nc{\Gtgr}{\Gt\Gr}
\nc{\Gt}{{\mathsf{G}_2}}
\nc{\OGr}{{\mathsf{OGr}}}
\nc{\Flag}{{\mathsf{Fl}}}
\nc{\Kosz}{{\mathsf{Kosz}}}
\nc{\LGr}{{\mathsf{LGr}}}
\nc{\LFl}{{\mathsf{LFl}}}
\nc{\SGr}{{\mathsf{SGr}}}
\nc{\OF}{{\mathsf{OF}}}
\nc{\Fl}{{\mathsf{Fl}}}
\nc{\Bl}{{\mathsf{Bl}}}

\nc{\Gm}{{\mathbb{G}_m}}
\nc{\GL}{{\mathsf{GL}}}
\nc{\PGL}{{\mathsf{PGL}}}
\nc{\GSL}{{\mathsf{SL}}}
\nc{\GSO}{{\mathsf{SO}}}
\nc{\SP}{{\mathsf{Sp}}}
\nc{\Spin}{{\mathsf{Spin}}}

\nc{\fsl}{{\mathfrak{sl}}}
\nc{\fso}{{\mathfrak{so}}}
\nc{\fgl}{{\mathfrak{gl}}}

\nc{\ev}{{\mathsf{ev}}}
\nc{\coev}{{\mathsf{coev}}}
\nc{\tr}{{\mathsf{tr}}}
\nc{\id}{{\mathsf{id}}}
\nc{\opp}{{\mathsf{opp}}}

\nc{\NHH}{{\mathop{\mathsf{NHH}}}}
\nc{\Cliff}{{\mathop{\mathcal{C}\!\ell}}}

\theoremstyle{plain}

\newtheorem{theorem}{Theorem}[section]
\newtheorem{conjecture}[theorem]{Conjecture}
\newtheorem{lemma}[theorem]{Lemma}
\newtheorem{proposition}[theorem]{Proposition}

\theoremstyle{definition}

\newtheorem{definition}[theorem]{Definition}

\newtheorem{example}[theorem]{Example}

\theoremstyle{remark}

\newtheorem{remark}[theorem]{Remark}
\newtheorem{question}[theorem]{Question}

\contact[akuznet@mi.ras.ru]{{\sloppy
\parbox{0.9\textwidth}{
Steklov Mathematical Institute,
8 Gubkin str., Moscow 119991 Russia
\\[5pt]
The Poncelet Laboratory, Independent University of Moscow
\hfill\\[5pt]
Laboratory of Algebraic Geometry, SU-HSE
}\bigskip}}

\title[Semiorthogonal decompositions in algebraic geometry]%
{Semiorthogonal decompositions\\in algebraic geometry}
\author{Alexander Kuznetsov}

\begin{document}

\maketitle

\begin{abstract}
In this review we discuss what is known about semiorthogonal decompositions 
of derived categories of algebraic varieties. We review existing constructions, 
especially the homological projective duality approach, and discuss some related
issues such as categorical resolutions of singularities.
\end{abstract}



\section*{Introduction}

In recent years an extensive investigation of semiorthogonal decompositions
of derived categories of coherent sheaves on algebraic varieties has been
done, and now we know quite a lot of examples and some general constructions.
With time it is becoming more and more clear that semiorthogonal components
of derived categories can be thought of as the main objects in noncommutative
algebraic geometry. In this paper I will try to review what is known in this direction ---
how one can construct semiorthogonal decompositions and how one can use them.

In section~\ref{sec-sod} we will recall the basic notions, 
discuss the most frequently used semiorthogonal decompositions,
and state the base change formula.
In section~\ref{sec-hpd} we review the theory of homological projective duality
which up to now is the most powerful method to construct semiorthogonal decompositions.
In section~\ref{sec-crs} we discuss categorical resolutions of singularities,
a subject interesting by itself, and at the same time inseparable from homological
projective duality.
In section~\ref{sec-ex} examples of homologically projectively dual varieties
are listed.
Finally, in section~\ref{sec-var} we discuss semiorthogonal decompositions
of varieties of small dimension. 

I should stress that in the area of algebraic geometry described in this paper there
are more questions than answers, but it really looks very promising.
Also, due to volume constraints I had to leave out
many interesting topics closely related to the main subject, such as 
the categorical Griffiths component, Hochschild homology and cohomology, and many others.

\section{Semiorthogonal decompositions}\label{sec-sod}

This paper can be considered as a continuation and a development
of the ICM 2002 talk \cite{BO02} of Alexei Bondal and Dmitri Orlov.
So I will freely use results and definitions from~\cite{BO02} and restrict
myself to a very short reminder of the most basic notion.
In particular, the reader is referred to~\cite{BO02} for the definition 
of a Serre functor, Fourier--Mukai transform, etc.

\subsection{A short reminder}

Recall that a {\sf semiorthogonal decomposition} of a triangulated category $\CT$ is
a collection $\CA_1,\dots,\CA_n$ of full triangulated subcategories 
such that:
$(1)$ for all $1 \le j < i \le n$
and any objects $A_i \in \CA_i$, $A_j \in \CA_j$ one has
$\Hom_\CT(A_i,A_j) = 0$;
$(2)$ the smallest triangulated subcategory of $\CT$ containing $\CA_1,\dots,\CA_n$
coincides with $\CT$.
We will use the notation $\CT = \langle \CA_1,\dots,\CA_n \rangle$ for a semiorthogonal 
decomposition of~$\CT$ with {\sf components} $\CA_1$, \dots, $\CA_n$.

We will be mostly interested in semiorthogonal decompositions of $\BD^b(\coh(X))$,
the {\em bounded derived category of coherent sheaves} on an algebraic variety $X$ 
which in most cases will be assumed to be smooth and projective over a base field $\kk$.

Recall that a full triangulated subcategory $\CA \subset \CT$ is {\sf admissible}
if its embedding functor $i:\CA \to \CT$ has both left and right adjoint functors 
$i^*,i^!:\CT \to \CA$. An admissible subcategory $\CA \subset \CT$ gives rise 
to a pair of semiorthogonal decompositions

\begin{equation}\label{lsod}
\CT = \langle \CA, {}^\perp\CA \rangle
\qquad\text{and}\qquad
\CT = \langle \CA^\perp, \CA \rangle,
\end{equation} 
where
\begin{align}
{}^\perp\CA := \{ T \in \CT\ |\ \Hom(T,A[t]) = 0\ \text{for all $A \in \CA$, $t \in \ZZ$} \},\label{lort}\\
\CA^\perp := \{ T \in \CT\ |\ \Hom(A[t],T) = 0\ \text{for all $A \in \CA$, $t \in \ZZ$} \},\label{rort}
\end{align} 
are the {\sf left} and the {\sf right orthogonals} to $\CA$ in $\CT$.
More generally, 
if $\CA_1,\dots,\CA_m$ is a semiorthogonal collection of admissible subcategories in $\CT$,
then for each $0 \le k \le m$ there is a semiorthogonal decomposition
\begin{equation}\label{gensod}
\CT = \langle \CA_1,\dots,\CA_k, 
{}^\perp\langle \CA_1,\dots,\CA_k \rangle \cap \langle \CA_{k+1}, \dots, \CA_m \rangle^\perp,
\CA_{k+1},\dots,\CA_m \rangle.
\end{equation}

The simplest example of an admissible subcategory is the one generated by an exceptional object.
Recall that an object $E$ is {\sf exceptional} if $\Hom(E,E) = \kk$ and $\Hom(E,E[t]) = 0$ 
for $t \ne 0$. An {\sf exceptional collection} is a collection of exceptional objects
$E_1,E_2,\dots,E_m$ such that $\Hom(E_i,E_j[t]) = 0$ for all $i > j$ and all $t \in \ZZ$.
An exceptional collection in $\CT$ gives rise to a semiorthogonal decomposition
\begin{equation}\label{sodec}
\CT = \langle \CA, E_1,\dots, E_m \rangle
\quad\text{with}\quad
\CA = \langle E_1,\dots, E_m \rangle^\perp.
\end{equation}
Here $E_i$ denotes the subcategory generated by the same named exceptional object.
If the category $\CA$ in~\eqref{sodec} is zero the exceptional collection is called {\sf full}.

\subsection{Full exceptional collections}

There are several well-known and quite useful semiorthogonal decompositions.
The simplest example is the following

\begin{theorem}[Beilinson's collection]
There is a full exceptional collection
\begin{equation}\label{beiexcol}
\BD^b(\coh(\PP^n)) = \langle \CO_{\PP^n}, \CO_{\PP^n}(1), \dots, \CO_{\PP^n}(n) \rangle.
\end{equation} 
\end{theorem}

Of course, twisting by $\CO_{\PP^n}(t)$ we get 
$\langle \CO_{\PP^n}(t),\CO_{\PP^n}(t+1),\dots,\CO_{\PP^n}(t+n) \rangle$
which is also a full exceptional collection for each $t \in \ZZ$.

A bit more general is the Grassmannian variety:

\begin{theorem}[Kapranov's collection, \cite{Kap92}]
Let $\Gr(k,n)$ be the Grassmannian of $k$-dimensional subspaces in a vector space of dimension $n$.
Let $\CU$ be the tautological subbundle of rank $k$. If $\chr\kk = 0$ then there is a semiorthogonal 
decomposition
\begin{equation}\label{kapexcol}
\BD^b(\coh(\Gr(k,n))) = \langle \ \Sigma^\alpha\CU^\vee\  \rangle_{\alpha \in R(k,n-k)},
\end{equation} 
where $R(k,n-k)$ is the $k\times(n-k)$ rectangle, $\alpha$ is a Young diagram, and $\Sigma^\alpha$
is the associated Schur functor.
\end{theorem}

When $\chr\kk > 0$ there is an exceptional collection as well, but it is a bit
more complicated, see~\cite{BLV}.

Another interesting case is the case of a smooth quadric $Q^n \subset \PP^{n+1}$.

\begin{theorem}[Kapranov's collection, \cite{Kap92}]
When $\chr\kk \ne 2$ there is a full exceptional collection
\begin{equation}\label{quexcol}
\BD^b(\coh(Q^n)) = 
\begin{cases}
\langle \SSS, \CO_{Q^n}, \CO_{Q^n}(1), \dots, \CO_{Q^n}(n-1) \rangle, & \text{ if $n$ is odd}\\
\langle \SSS^-, \SSS^+, \CO_{Q^n}, \CO_{Q^n}(1), \dots, \CO_{Q^n}(n-1) \rangle, & \text{ if $n$ is even}
\end{cases}
\end{equation}
where $\SSS$ and $\SSS^\pm$ are the spinor bundles.
\end{theorem}

Many exceptional collections have been constructed on other rational homogeneous spaces,
see e.g. \cite{S01}, \cite{K08b}, \cite{PS}, \cite{M}, \cite{FM}, and \cite{KP11}.
Full exceptional collections on smooth toric varieties (and stacks) were constructed
by Kawamata~\cite{Kaw}. Also exceptional collections were constructed on del Pezzo surfaces \cite{O92},
some Fano threefolds~\cite{O91,K96} and many other varieties. 

\subsection{Relative versions}

Let $S$ be a scheme and $E$ a vector bundle of rank~$r$ on it. Let $\PP_S(E)$
be its projectivization, $f:\PP_S(E) \to S$ the projection, and $\CO_{\PP_S(E)/S}(1)$
the Grothendieck line bundle on $\PP_S(E)$.

\begin{theorem}[\cite{O92}]\label{thmorlsod}
For each $i \in \ZZ$ the functor 
\begin{equation}\label{phii}
\Phi_i:\BD^b(\coh(S)) \to \BD^b(\coh(\PP_S(E))), 
\qquad
F \mapsto Lf^*(F)\lotimes\CO_{\PP_S(E)/S}(i)
\end{equation}
is fully faithful, and there is a semiorthogonal decomposition
\begin{equation}\label{orlsod}
\BD^b(\coh(\PP_S(E))) = 
\langle \Phi_0(\BD^b(\coh(S))),
\dots,\Phi_{r-1}(\BD^b(\coh(S))) \rangle.
\end{equation} 
\end{theorem}

Of course, analogously to the case of a projective space, one can replace the sequence 
of functors $\Phi_0,\dots,\Phi_{r-1}$ by $\Phi_t,\dots,\Phi_{t+r-1}$ for any $t \in \ZZ$.

An interesting new feature appears for Severi--Brauer varieties. Recall that a Severi--Brauer
variety over $S$ is a morphism $f:X \to S$ which \'etale locally is isomorphic to a projectivization
of a vector bundle. A Severi--Brauer variety $X$ can be constructed from a torsion element 
in the Brauer group $\Br(S)$ of $S$.

\begin{theorem}[Bernardara's decomposition, \cite{Be}]\label{thmbersod}
Let $f:X \to S$ be a Severi--Brauer variety of relative dimension $n$ 
and $\beta \in \Br(S)$ its Brauer class.
Then for each $i \in \ZZ$ there is a fully faithful functor 
$\Phi_i:\BD^b(\coh(S,\beta^i)) \to \BD^b(\coh(X))$ and 
a semiorthogonal decomposition
\begin{equation}\label{berexcol}
\hspace{-.8em}
\BD^b(\coh(X)) = \langle \Phi_0(\BD^b(\coh(S))),
\Phi_1(\BD^b(\coh(S,\beta))),
\dots,
\Phi_n(\BD^b(\coh(S,\beta^n)))\rangle.
\end{equation} 
\end{theorem}

Here $\coh(S,\beta^i)$ is the category of $\beta^i$-twisted coherent sheaves on $S$
and the functor $\Phi_i$ is given by $F \mapsto Lf^*(F)\lotimes\CO_{X/S}(i)$,
where the sheaf $\CO_{X/S}(i)$ is well defined as a $f^*\beta^{-i}$-twisted sheaf.

Another important semiorthogonal decomposition can be constructed for a smooth blowup.
Let $Y \subset X$ be a locally complete intersection subscheme of codimension $c$
and let $\tilde{X}$ be the blowup of $X$ with center in $Y$. Let $f:\tilde{X} \to X$
be the blowup morphism and $D \subset \tilde{X}$ the exceptional divisor of the blowup.
Let $i:D \to \tilde{X}$ be the embedding and $p:D \to Y$ the natural projection
(the restriction of $f$ to $D$). Note that $D \cong \PP_Y(\CN_{Y/X})$ is 
the projectivization of the normal bundle.

\begin{theorem}[Blowup formula, \cite{O92}]\label{thmbupsod}
The functor $Lf^*:\BD^b(\coh(X)) \to \BD^b(\coh(\tilde{X}))$
as well as the functors
\begin{equation*}
\Psi_k:\BD^b(\coh(Y)) \to \BD^b(\coh(\tilde{X})), 
\qquad
F \mapsto Ri_*(Lp^*(F)\lotimes\CO_{D/Y}(k)),
\end{equation*}
are fully faithful for all $k \in \ZZ$, and 
there is a semiorthogonal decomposition
\begin{equation}\label{blupsod}
\hspace{-.27em}
\BD^b(\coh(\tilde{X})) = 
\langle 
Lf^*(\BD^b(\coh(X))),
\Psi_0(\BD^b(\coh(Y))),
\dots,\Psi_{c-2}(\BD^b(\coh(Y))) \rangle.
\end{equation} 
\end{theorem}

Finally, consider a flat fibration in quadrics $f:X \to S$.
In other words, assume that $X \subset \PP_S(E)$ is a divisor 
of relative degree $2$ in a projectivization of a vector bundle $E$ of rank $n+2$
on a scheme $S$ corresponding to a line subbundle $\CL \subset S^2E^\vee$.

\begin{theorem}[Quadratic fibration formula, \cite{K08a}]\label{thm-qu-sod}
For each $i \in \ZZ$ there is a fully faithful functor 
\begin{equation*}
\Phi_i:\BD^b(\coh(S)) \to \BD^b(\coh(X)),
\qquad
F \mapsto Lf^*(F) \lotimes \CO_{X/S}(i)
\end{equation*}
and a semiorthogonal decomposition
\begin{equation}
\BD^b(\coh(X)) = 
\langle \BD^b(\coh(S,\Cliff_0)), \Phi_0(\BD^b(\coh(S))),\dots,\Phi_{n-1}(\BD^b(\coh(S))) \rangle,
\end{equation} 
where $\Cliff_0$ is the sheaf of even parts of Clifford algebras on $S$ associated with 
the quadric fibration $X \to S$.
\end{theorem}

The sheaf $\Cliff_0$ is a sheaf of $\CO_S$-algebras which as an $\CO_S$-module is isomorphic to
\begin{equation*}
\Cliff_0 \cong \CO_S \oplus (\Lambda^2E\otimes\CL) \oplus (\Lambda^4E\otimes\CL^2) \oplus \dots
\end{equation*}
and equipped with an algebra structure via the Clifford multiplication. If the dimension $n$
of fibers of $X \to S$ is odd, then $\Cliff_0$ is a sheaf of Azumaya algebras on the open subset of $S$
corresponding to nondegenerate quadrics (which of course may be empty). On the other hand, if
$n$ is even then $\CO_S \oplus \Lambda^nE \otimes \CL^{n/2}$ is a central subalgebra in $\Cliff_0$,
so the latter gives a sheaf $\tilde{\Cliff_0}$ of algebras on the twofold covering 
\begin{equation}\label{stilde}
\tilde{S}:= \Spec_S(\CO_S \oplus \Lambda^nE \otimes \CL^{n/2}) 
\end{equation} of $S$, and $\tilde{\Cliff_0}$
is a sheaf of Azumaya algebras on the preimage of the open subset of $S$ corresponding to
nondegenerate quadrics.

\subsection{Base change}

A triangulated category $\CT$ is {\sf $S$-linear} if it is equipped
with a module structure over the tensor triangulated category $\BD^b(\coh(S))$.
In particular, if $X$ is a scheme over $S$ and $f:X \to S$ is the structure morphism then
a semiorthogonal decomposition
\begin{equation}\label{sodx}
\BD^b(\coh(X)) = \langle \CA_1, \dots, \CA_m \rangle 
\end{equation}
is {\sf $S$-linear} if each of the subcategories $\CA_k$ is closed under tensoring 
with an object of $\BD^b(\coh(S))$, i.e.\ for $A \in \CA_k$ and $F \in \BD^b(\coh(S))$
one has $A \lotimes Lf^*(F) \in \CA_k$.

The semiorthogonal decompositions 
of Theorems~\ref{thmorlsod}, \ref{thmbersod} and~\ref{thm-qu-sod} are $S$-linear,
and the blowup formula of Theorem \ref{thmbupsod} is $X$-linear.
The advantage of linear semiorthogonal decompositions lies in the fact that they obey a 
base change result. For a base change $T \to S$ denote by $\pi:X\times_S T \to X$
the induced projection.

\begin{theorem}[\cite{K11}]\label{basethm}
If $X$ is an algebraic variety over $S$ and~\eqref{sodx}
is an $S$-linear semiorthogonal decomposition then for a change of base morphism $T \to S$
there is, under a certain technical condition, a $T$-linear semiorthogonal decomposition
\begin{equation*}
\BD^b(\coh(X\times_S T)) = \langle \CA_{1T}, \dots, \CA_{mT} \rangle
\end{equation*}
such that $\pi^*(A) \subset \CA_{iT}$ for any $A \in \CA_i$ and 
$\pi_*(A') \subset \CA_i$ for any $A' \in \CA_{iT}$ which has proper support over $X$.
\end{theorem}

\subsection{Important questions}

There are several questions which might be crucial for further investigations.

\begin{question}
Find a good condition for an exceptional collection to be full.
\end{question}

One might hope that if the collection generates the Grothendieck group
(or the Hochschild homology) of the category then it is full. However,
recent examples of quasiphantom and phantom categories (see section~\ref{ssec-surfaces}) 
show that this is not the case. Still we may hope that in the categories generated by 
exceptional collections there are no phantoms. In other words
one could hope that the following is true.

\begin{conjecture}
Let $\CT = \langle E_1,\dots, E_n \rangle$ be a triangulated category generated 
by an exceptional collection. Then any exceptional collection of length $n$ in $\CT$ is full.
\end{conjecture}

If there is an action of a group $G$ on an algebraic variety $X$, one can consider
the equivariant derived category $\BD^b(\coh^G(X))$ along with the usual derived 
category. In many interesting cases (flag varieties, toric varieties, GIT quotients)
it is quite easy to construct a full exceptional collection in the equivariant category.
It would be extremely useful to find a way to transform it into a full exceptional
collection in the usual category. In some sense the results of~\cite{KP11}
can be considered as an example of such an approach.

Another very important question is to find possible restrictions for existence
of semiorthogonal decompositions. Up to now there are only several cases when we can {\em prove}
indecomposability of a category. The first is the derived category of a curve of positive genus.
The proof (see e.g.~\cite{Ok11}) is based on special properties of categories of homological 
dimension~1. Another is the derived category of a Calabi--Yau variety (smooth connected
variety with trivial canonical class). Its indecomposability is proved by a surprisingly simple 
argument due to Bridgeland~\cite{Br99}. This was further generalized in~\cite{KaOk12}
to varieties with globally generated canonical class. On the other hand, 
the original argument of Bridgeland generalizes to any connected Calabi--Yau category
(i.e.~with the Serre functor isomorphic to a shift and Hochschild cohomology in degree 
zero isomorphic to $\kk$).

\section{Homological projective duality}\label{sec-hpd}

The starting point of a homological projective duality (HP duality for short) is a smooth
projective variety $X$ with a morphism into a projective space and a semiorthogonal 
decomposition of $\BD^b(\coh(X))$ of a very special type.

\subsection{Lefschetz decompositions}

Let $X$ be an algebraic variety and $\CL$ a line bundle on $X$.

\begin{definition}\label{defld}
A {\sf right Lefschetz decomposition} of $\BD^b(\coh(X))$ with respect to $\CL$
is a semiorthogonal decomposition of form
\begin{equation}\label{rld}
\BD^b(\coh(X)) = \langle \CA_0, \CA_1\otimes \CL, \dots, \CA_{m-1}\otimes \CL^{m-1} \rangle
\end{equation}
with
$0 \subset \CA_{m-1} \subset \dots \subset \CA_1 \subset \CA_0$.
In other words, each component of the decomposition is a subcategory of the previous component 
twisted by $\CL$.

Analogously, a {\sf left Lefschetz decomposition} of $\BD^b(\coh(X))$ with respect to $\CL$
is a semiorthogonal decomposition of form
\begin{equation}\label{lld}
\BD^b(\coh(X)) = \langle \CB_{m-1}\otimes\CL^{1-m}, \dots, \CB_1\otimes\CL^{-1}, \CB_0 \rangle
\end{equation}
with
$0 \subset \CB_{m-1} \subset \dots \subset \CB_1 \subset \CB_0$.
\end{definition}

The subcategories $\CA_i$ (resp.\ $\CB_i$) forming a Lefschetz decomposition will be called {\sf blocks},
the largest will be called {\sf the first block}. Usually we will consider right Lefschetz decompositions. 
So, we will call them simply Lefschetz decompositions.

Beilinson's collection on $\PP^n$ is an example of a Lefschetz decomposition with
$\CA_0 = \CA_1 = \dots = \CA_n = \langle \CO_{\PP^n} \rangle$. Kapranov's collection
on the Grassmannian $\Gr(k,n)$ also has a Lefschetz structure with the category $\CA_i$ 
generated by $\Sigma^\alpha\CU^\vee$ for $\alpha \in R(k-1,n-k-i)$.

Note that in Definition~\ref{defld} one can replace the twist by a line bundle 
with any other autoequivalence of $\BD^b(\coh(X))$ and get the notion of a Lefschetz 
decomposition with respect to an autoequivalence. This may be especially useful 
when dealing with arbitrary triangulated categories.

It is also useful to know that for a given line bundle $\CL$ a Lefschetz decomposition 
is completely determined by its first block. Moreover, an admissible subcategory 
extends to a right Lefschetz decomposition if and only if it extends to a left Lefschetz decomposition.
The simplest example of an admissible subcategory which does not extend to a Lefschetz decomposition
is the subcategory $\langle \CO_{\PP^2}, \CO_{\PP^2}(2) \rangle \subset \BD^b(\coh(\PP^2))$.

\begin{question}
Find a good sufficient condition for a Lefschetz extendability 
of an admissible subcategory $\CA_0 \subset \BD^b(\coh(X))$.
\end{question}

One can define a partial ordering on the set of all Lefschetz decompositions of $\BD^b(\coh(X))$ 
by inclusions of their first blocks. As we will see soon, the most interesting and strong results 
are obtained by using {\sf minimal} Lefschetz decompositions.

\subsection{Hyperplane sections}

Let $X$ be a smooth projective variety with a morphism into a projective space $f:X \to \PP(V)$ 
(not necessarily an embedding). Put $\CO_X(1) := f^*\CO_{\PP(V)}(1)$ and assume that a right
Lefschetz decomposition with respect to $\CO_X(1)$ 
\begin{equation}\label{ld}
\BD^b(\coh(X)) = \langle \CA_0, \CA_1(1), \dots, \CA_{m-1}(m-1) \rangle
\end{equation}
is given (we abbreviate $\CA_i(i) := \CA_i \otimes\CO_X(i)$). 
Consider the dual projective space $\PP(V^\vee)$. By the base change (Theorem~\ref{basethm})
the product $X \times \PP(V^\vee)$ inherits a $\PP(V^\vee)$-linear semiorthogonal decomposition
\begin{equation*}
\BD^b(\coh(X\times\PP(V^\vee))) = 
\langle \CA_{0\PP(V^\vee)}, \CA_1(1)_{\PP(V^\vee)}, \dots, \CA_{m-1}(m-1)_{\PP(V^\vee)} \rangle
\end{equation*}
Consider the {\sf universal hyperplane section of $X$},
$\CX := X \times_{\PP(V)} Q \subset X \times \PP(V^\vee)$,
where $Q \subset \PP(V)\times\PP(V^\vee)$ is the {\sf incidence quadric} and
denote by $\alpha:\CX \to X\times\PP(V^\vee)$ the natural embedding. 

\begin{lemma}
The functor $L\alpha^*:\BD^b(\coh(X\times\PP(V^\vee))) \to \BD^b(\coh(\CX))$
is fully faithful on each of the subcategories $\CA_1(1)_{\PP(V^\vee)}$, \dots,
$\CA_{m-1}(m-1)_{\PP(V^\vee)}$ and induces a $\PP(V^\vee)$-linear 
semiorthogonal decomposition
\begin{equation}\label{hpdcat}
\BD^b(\coh(\CX)) = \langle \CC,  \CA_1(1)_{\PP(V^\vee)}, \dots, \CA_{m-1}(m-1)_{\PP(V^\vee)} \rangle.
\end{equation}
\end{lemma}

The first component $\CC$ of this decomposition is called the {\sf HP dual category} of~$X$. 
It is a very interesting category, especially if it can be identified with the derived 
category of some algebraic variety~$Y$. In this case this variety is called the 
{\sf HP dual variety} of $X$.

\begin{definition}\label{defhpd}
An algebraic variety $Y$ equipped with a morphism $g:Y \to \PP(V^\vee)$ is called
{\sf homologically projectively dual} to $f:X \to \PP(V)$ with respect to a given
Lefschetz decomposition~\eqref{ld} if there is given an object 
$\CE \in \BD^b(\coh(Q(X,Y)))$ such that the Fourier--Mukai functor
$\Phi_\CE:\BD^b(\coh(Y)) \to \BD^b(\coh(\CX))$ is an equivalence 
onto the HP dual subcategory $\CC \subset \BD^b(\coh(\CX))$ of~\eqref{hpdcat}.
\end{definition}

Here
$Q(X,Y) = (X\times Y) \times_{\PP(V)\times\PP(V^\vee)} Q = \CX \times_{\PP(V^\vee)} Y$.
If a homological projective duality between varieties $X$ and $Y$ is established then there
is an interesting relation between derived categories of their linear sections.

\subsection{HP duality statement}

For each linear subspace $L \subset V^\vee$ denote by
$L^\perp := \Ker(V \to L^\vee)$
its orthogonal complement in $V$. Further denote 
\begin{equation}\label{defxlyl}
X_L := X \times_{\PP(V)}\PP(L^\perp),
\qquad
Y_L := Y \times_{\PP(V^\vee)} \PP(L).
\end{equation}
Varieties defined in this way are called {\sf mutually orthogonal linear sections} of $X$ and~$Y$.
We will say that $X_L$ and $Y_L$ {\sf have expected dimensions} if
\begin{equation*}
\dim X_L = \dim X - r
\qquad\text{and}\qquad
\dim Y_L = \dim Y - (N - r),
\end{equation*}
where $N = \dim V$ and $r = \dim L$ (so that $N-r = \dim L^\perp$).

\begin{theorem}[Homological projective duality, \cite{K07a}]\label{thmhpd}
Let $(Y,g)$ be an HP dual variety for $(X,f)$ with respect to~\eqref{ld}. Then
\begin{enumerate}
\item 
$Y$ is smooth and $\BD^b(\coh(Y))$ has an admissible subcategory $\CB_0$
equivalent to $\CA_0$ and extending to a left Lefschetz decomposition
\begin{equation}\label{lddual}
\hspace{-2em}
\BD^b(\coh(Y)) = \langle \CB_{n-1}(1-n), \dots, \CB_1(-1), \CB_0 \rangle,
\quad
\CB_{n-1} \subset \dots \subset \CB_1 \subset \CB_0.
\end{equation} 
\item 
$(X,f)$ is HP dual to $(Y,g)$ with respect to~\eqref{lddual}.
\item 
The set of critical values of $g$ is the classical projective dual of $X$.
\item 
For any subspace $L \subset V^\vee$ if $X_L$ and $Y_L$ have expected dimensions 
then there are semiorthogonal decompositions
\begin{align}
\BD^b(\coh(X_L)) &= \langle \CC_L, \CA_r(r), \dots, \CA_{m-1}(m-1) \rangle,\label{dxl}\\
\BD^b(\coh(Y_L)) &= \langle \CB_{n-1}(1-n), \dots, \CB_{N-r}(r-N), \CC_L \rangle\label{dyl}
\end{align} 
with the same triangulated category $\CC_L$ appearing in the RHS.
\end{enumerate}
\end{theorem}

The decomposition~\eqref{lddual} of $\BD^b(\coh(Y))$ will be referred to 
as the {\sf HP dual Lefschetz decomposition}. The common component $\CC_L$ 
of decompositions~\eqref{dxl} and~\eqref{dyl}
will be referred to as the {\sf nontrivial part} of the derived categories 
of $X_L$ and $Y_L$, while the subcategories $\CA_i(i)$ and $\CB_j(-j)$ 
(one checks that the pullback functors for the embeddings $X_L \to X$ and $Y_L \to Y$
are fully faithful on the subcategories $\CA_i$ and $\CB_j$ for $i \ge r$ and $j \ge N-r$)
are considered as trivial (in the sense that they come from the ambient varieties).

The first two statements of this Theorem show that the relation we are dealing with
is indeed a duality, the third statement shows the relation to the classical projective 
duality (and so justifies the word ``projective'' in the name), and the last statement
is the real result. We will soon see how powerful it is.

Note also that in the statement of the Theorem the linear sections $X_L$ and $Y_L$ need not be smooth.
In fact, one can show that for HP dual varieties $X$ and $Y$ a section $X_L$ is smooth if and only if 
its orthogonal section $Y_L$ is smooth, but no matter whether this is the case or not, 
the decompositions~\eqref{dxl} and~\eqref{dyl} hold true.

Now let us say some words about the relations of the Lefschetz decompositions~\eqref{ld}
and~\eqref{lddual} for HP dual varieties. As it was already mentioned, the largest components of those
are just equivalent
$\CB_0 \cong \CA_0$.
Further, the component $\CB_i$ is very close to the orthogonal complement of $\CA_{N-1-i}$ in~$\CA_0$.
More precisely, these two categories have semiorthogonal decompositions with the same components
but with in general different gluing functors. This can be visualized by 
a picture.
\begin{align*}
\ytableausetup{nosmalltableaux}
\begin{ytableau}
\none[\textstyle \CA_{0}] & \none[\textstyle \CA_{1}] & 
\none[\textstyle \CA_{2}] & \none[\textstyle \CA_{3}] & 
\none[\textstyle \CA_{4}] & \none[\textstyle \CA_{5}] & 
\none[\textstyle \CA_{6}] & \none[\textstyle \CA_{7}] & 
\none[\textstyle \CA_{8}] & \none[\textstyle \CA_{9}] \\
*(lightgray)&*(lightgray)&*(lightgray)&*(lightgray)&
*(lightgray)&*(lightgray)&*(lightgray)&*(lightgray)&
*(lightgray)&*(lightgray)&&&&&&&&&&\\
*(lightgray)&*(lightgray)&*(lightgray)&*(lightgray)&
*(lightgray)&*(lightgray)&*(lightgray)&*(lightgray)&&&&&&&&&&&&\\
*(lightgray)&*(lightgray)&*(lightgray)&*(lightgray)&
*(lightgray)&&&&&&&&&&&&&&&\\
*(lightgray)&*(lightgray)&*(lightgray)&&&&&&&&&&&&&&&&&\\
\none & \none & \none & 
\none[\textstyle \CB_{16}] & \none[\textstyle \CB_{15}] & 
\none[\textstyle \CB_{14}] & \none[\textstyle \CB_{13}] & 
\none[\textstyle \CB_{12}] & \none[\textstyle \CB_{11}] & 
\none[\textstyle \CB_{10}] & \none[\textstyle \CB_{9}] & 
\none[\textstyle \CB_{8}] & \none[\textstyle \CB_{7}] & 
\none[\textstyle \CB_{6}] & \none[\textstyle \CB_{5}] & 
\none[\textstyle \CB_{4}] & \none[\textstyle \CB_{3}] & 
\none[\textstyle \CB_{2}] & \none[\textstyle \CB_{1}] & 
\none[\textstyle \CB_{0}] 
\end{ytableau}
\end{align*}
The gray part of the picture corresponds to the initial Lefschetz decomposition, 
the columns correspond to its blocks, while the white part corresponds to the dual
decomposition, the complementary columns correspond to the complementary subcategories
of the dual Lefschetz decomposition. The number of rows is equal to the number of different 
components in the initial (and the dual) Lefschetz decomposition. In this example picture 
$\CA_0 = \CA_1 = \CA_2 \ne \CA_3 = \CA_4 \ne \CA_5 = \CA_6 = \CA_7 \ne \CA_8 = \CA_9$,
and so one can say that the rows correspond to the ``primitive parts''
$(\CA_3)^\perp_{\CA_0}$, $(\CA_5)^\perp_{\CA_3}$, $(\CA_8)^\perp_{\CA_5}$, and $\CA_8$
of all the categories in the picture,
the length of the initial decomposition is $m = 10$, the length of the dual decomposition 
is $n = 17$, while the dimension of the ambient space is $N = 20$.

Note that $\CB_i = 0$ if and only 
if $\CA_{N-1-i} = \CA_0$, so the number $n$ of components in~\eqref{lddual} equals $N$ minus
the number of components in~\eqref{ld} equal to $\CA_0$.

In fact, the best (in many aspects) situation is when in the original Lefschetz 
decomposition~\eqref{ld} all components coincide $\CA_0 = \CA_1 = \dots = \CA_{m-1}$
(such Lefschetz decompositions are called {\sf rectangular}). 
Then the HP dual Lefschetz decomposition is also rectangular, has the same components 
$\CB_0 = \CB_1 = \dots = \CB_{n-1} \cong \CA_0$ and
\begin{equation*}
n = N - m
\end{equation*}
(in particular in a picture analogous to the above the gray and the white parts
are rectangles, which explains the name ``rectangular'').
Moreover, in this case for any $0 < r < N$ one has either $r \ge m$ or $N-r \ge n$,
hence in decompositions~\eqref{dxl} and~\eqref{dyl} either the first or the second
category has only one component $\CC_L$ and nothing else. Then the other
decomposition shows that the nontrivial component of the derived category of 
a linear section is equivalent to the derived category of the orthogonal linear section 
of the dual variety.

\subsection{HP duality and noncommutative varieties}\label{ss-hpd-noncom}

In general, the HP dual category
$\CC \subset \BD^b(\coh(\CX))$ defined by~\eqref{hpdcat} need not be equivalent to $\BD^b(\coh(Y))$
for an algebraic variety $Y$. In fact, only a few such cases are known ---
the linear duality, the duality for quadrics, the duality for Grassmannians $\Gr(2,4)$ and $\Gr(2,5)$,
and the spinor variety $\SSSS_5$ (see section~\ref{sec-ex}).

One can get many additional interesting examples by allowing $Y$ to be a noncommutative variety.
Here a noncommutative variety can be understood in different ways. If one uses the most general 
sense --- as a semiorthogonal component of the derived category of an algebraic variety ---
then tautologically the HP dual category $\CC$ itself will provide a noncommutative HP dual variety.
In fact, one can develop a theory of HP duality for noncommutative varieties in this sense 
and prove the same results (see~\cite{K07b}). However, in this
most general form the semiorthogonal decompositions provided by the HP duality Theorem
will not have an apparent geometric interpretation. 

In fact, an interesting geometry arises in HP duality if the dual variety $Y$
is close to a commutative variety. For example, it often happens that there is a (commutative) 
algebraic variety $Y_0$ with a map $g_0:Y_0 \to \PP(V^\vee)$, a sheaf of 
finite $\CO_{Y_0}$-algebras $\CR$ on $Y_0$ whose bounded derived
category $\BD^b(\coh(Y_0,\CR))$ of coherent $\CR$-modules on $Y_0$ is equivalent 
to the HP dual category $\CC$ of $X$
and such that the equivalence $\CC \cong \BD^b(\coh(Y_0,\CR))$ is given by an appropriate object
$\CE \in \BD^b(\coh(Q(X,Y_0),\CR))$. Of course, one can easily allow here the variety $X$
also to be noncommutative in the same sense. It is easy to modify all the definitions accordingly.

In section~\ref{sec-ex} we discuss examples showing that this generalization is meaningful.
Among such examples are the Veronese--Clifford duality, the Grassmannian--Pfaffian duality,
and their generalizations.

In fact, in some of these examples, the HP duality Theorem~\ref{thmhpd} still gives 
semiorthogonal decompositions for usual commutative varieties (even though the dual 
variety is noncommutative). Indeed, the sheaf of algebras $\CR$ on $Y_0$ is frequently isomorphic 
to a matrix algebra on an open subset of $Y_0$, typically, on its smooth locus ---
in fact, in these cases the noncommutative variety $(Y_0,\CR)$ can be thought of
as a categorical resolution of singularities of $Y$. In this situation,
taking a subspace $L \subset V^\vee$ such that $Y_{0L}$ is contained in that 
open subset, one gets $\BD^b(\coh(Y_L)) = \BD^b(\coh(Y_{0L},\CR)) \cong \BD^b(\coh(Y_{0L}))$.

\section{Categorical resolutions of singularities}\label{sec-crs}

As it was explained above (and we will see in some of the examples below) 
in many cases the HP dual variety looks as a noncommutative (or categorical) 
resolution of singularities of a singular variety. So, a good notion 
of a categorical resolution is necessary for the theory.

\subsection{The definition}

If $\pi:\TY \to Y$ is a resolution of singularities, we have
an adjoint pair of triangulated functors $R\pi_*:\BD^b(\coh(\TY)) \to \BD^b(\coh(Y))$ and
$L\pi^*:\BD^\perf(Y) \to \BD^b(\coh(\TY))$ (here $\BD^\perf(Y)$ stands for the category
of perfect complexes on $Y$). We axiomatize this situation in the following

\begin{definition}[cf.\ \cite{K08c,KL12}]\label{defcrs}
A {\sf categorical resolution of singularities} of a scheme $Y$ is 
a smooth triangulated category $\CT$ and an adjoint pair 
of triangulated functors $\pi_*:\CT \to \BD^b(\coh(Y))$ and 
$\pi^*:\BD^\perf(Y) \to \CT$ 
such that
$\pi_*\circ\pi^* \cong \id_{\BD^\perf(Y)}$. In particular,
the functor $\pi^*$ is fully faithful.
\end{definition}

We will not discuss the notion of smoothness for a triangulated category.
In fact, for our purposes it is always enough to assume that $\CT$ is an admissible 
$Y$-linear subcategory of $\BD^b(\coh(\TY))$ for a geometric resolution $\TY \to Y$.

Let $(\CT,\pi_*,\pi^*)$ and $(\CT',\pi'_*,{\pi'}^*)$ be two categorical resolutions of $Y$.
We will say that $\CT$ {\sf dominates} $\CT'$ if there is a fully faithful functor
$\epsilon:\CT' \to \CT$ such that $\pi'_* = \pi_*\circ\epsilon$.
Clearly, this is compatible with the usual dominance relation between 
geometric resolutions --- if a resolution $\pi:\TY \to Y$ factors as 
$\TY \xrightarrow{\ f\ } \TY' \xrightarrow{\ \pi'\ } Y$ then the functor 
$\epsilon := Lf^*:\BD^b(\coh(\TY')) \to \BD^b(\coh(\TY))$ is fully faithful and
\begin{equation*}
R\pi_*\circ Lf^* = R\pi'_*\circ Rf_*\circ Lf^* \cong R\pi'_*.
\end{equation*}

Categorical resolutions have two advantages in comparison with geometric ones.
First, if $Y$ has {\sf irrational singularities} the pullback functor 
for a geometric resolution is {never} fully faithful and so its derived
category is not a categorical resolution in sense of Definition~\ref{defcrs}. 
However, it was shown in~\cite{KL12} that any separated scheme of finite type 
(even nonreduced) over a field of zero characteristic admits a categorical resolution.

The second advantage is that the dominance order for categorical resolutions
is more flexible. For example, in many examples one can find a categorical
resolution which is much smaller than any geometric resolution.
There are strong indications that the Minimal Model Program on the categorical level
may be much simpler than the classical one. In particular, we expect the following.

\begin{conjecture}[cf.\ \cite{BO02}]\label{mincr}
For any quasiprojective scheme $Y$ there exists a categorical resolution 
which is minimal with respect to the dominance order.
\end{conjecture}

\subsection{Examples of categorical resolutions}

As an evidence for the conjecture we will construct 
categorical resolutions which are presumably minimal.

\begin{theorem}[\cite{K08c}]\label{ldtocr}
Let $f:\TY \to Y$ be a resolution of singularities and let $E$ be
the exceptional divisor with $i:E \to \TY$ being the embedding. 
Assume that the derived category $\BD^b(\coh(E))$ has a left Lefschetz decomposition 
with respect to the conormal bundle $\CO_E(-E)$:
\begin{equation}\label{lde}
\BD^b(\coh(E)) = \langle \CC_{m-1}((m-1)E), \dots, \CC_1(E), \CC_0 \rangle,
\end{equation}
which is $Y$-linear and has 
$Li^*(Lf^*(\BD^\perf(Y))) \subset \CC_0$.
Then the functor $Ri_*$ is fully faithful on subcategories $\CC_k \subset \BD^b(\coh(E))$
for $k > 0$, the subcategory 
\begin{equation}\label{cres}
\TCC := \{ F \in \BD^b(\coh(\TY))\ |\ Li^*(F) \in \CC_0 \}
\end{equation}
is admissible in $\BD^b(\coh(\TY))$, and there is a semiorthogonal decomposition
\begin{equation}\label{sodty}
\BD^b(\coh(\TY)) = \langle Ri_*(\CC_{m-1}((m-1)E)), \dots, Ri_*(\CC_1(E)), \TCC \rangle.
\end{equation}
Moreover, the functor $Lf^*:\BD^\perf(Y) \to \BD^b(\coh(\TY))$ factors as a composition
of a fully faithful functor $\pi^*:\BD^\perf(Y) \to \TCC$ with the embedding 
$\gamma:\TCC \to \BD^b(\coh(\TY))$,
and the functors $\pi_* := Rf_*\circ\gamma$ and $\pi^*$
give $\TCC$ a structure of a categorical resolution of singularities of $Y$.
\end{theorem}

If $\CC'_0 \subset \CC_0 \subset \BD^b(\coh(E))$ are two admissible Lefschetz extendable 
subcategories (with respect to the conormal bundle) then clearly by~\eqref{cres} 
the categorical resolution $\TCC'$ constructed from $\CC'_0$ is a subcategory 
in the categorical resolution $\TCC$ constructed from $\CC_0$. Moreover,
if $\epsilon:\TCC' \to \TCC$ is the embedding functor then $\pi'_* = \pi_*\circ\epsilon$,
so $\TCC$ dominates $\TCC'$. This shows that minimal categorical resolutions
are related to minimal Lefschetz decompositions.

As an example of the application of the above Theorem consider the cone $Y$ over a smooth
projective variety $X \subset \PP(V)$. Then $\TY = \Tot_X(\CO_X(-1))$, the total
space of the line bundle $\CO_X(-1) = \CO_{\PP(V)}(-1)_{|X}$, 
is a geometric resolution of $Y$.
The exceptional divisor of the natural morphism $f:\TY \to Y$ then identifies with 
the zero section of the total space, $E = X$, and the conormal bundle identifies
with~$\CO_X(1)$. So, a left Lefschetz decomposition of $\BD^b(\coh(X))$ with respect
to~$\CO_X(1)$ gives a categorical resolution of the cone $Y$ over $X$.

\begin{example}\label{vercone}
Take $X = \PP^3$ with the double Veronese embedding $f:\PP^3 \to \PP^9$,
so that $f^*\CO_{\PP^9}(1) = \CO_{\PP^3}(2)$, and a left Lefschetz decomposition 
\begin{equation*}
\BD^b(\coh(\PP^3)) = \langle \CC_1(-2), \CC_0 \rangle
\quad\text{with $\CC_0 = \CC_1 = \langle \CO_{\PP^3}(-1), \CO_{\PP^3} \rangle$.}
\end{equation*}
Then the category 
$\TCC := \{ F \in \BD^b(\coh(\Tot_{\PP^3}(\CO_{\PP^3}(-2))))\ |\ 
Li^*F \in \langle \CO_{\PP^3}(-1), \CO_{\PP^3} \rangle \}$
is a categorical resolution of the Veronese cone, which is significantly smaller 
than the usual geometric resolution. It is expected to be minimal.
\end{example}

\subsection{Crepancy of categorical resolutions}

Crepancy is an important property of a resolution which in the geometric
situation ensures its minimality. A resolution $f:\TY \to Y$
is {\sf crepant} if the relative canonical class $K_{\TY/Y}$ is trivial.
There is an analogue of crepancy for categorical resolutions. 
In fact, there are two such analogues.

\begin{definition}[\cite{K08c}]
A categorical resolution $(\CT,\pi_*,\pi^*)$ of a scheme $Y$ is {\sf weakly crepant}
if the functor $\pi^*:\BD^\perf(Y) \to \CT$ is both left and right adjoint 
to the functor $\pi_*:\CT \to \BD^b(\coh(Y))$.
\end{definition}

By Grothendieck duality, the right adjoint of the derived pushforward functor 
$Rf_*:\BD^b(\coh(\TY)) \to \BD^b(\coh(Y))$ is given by $f^!(F) = Lf^*(F)\otimes \CO_\TY(K_{\TY/Y})$,
so for a geometric resolution crepancy and weak crepancy are equivalent.

\begin{definition}[\cite{K08c}]
A categorical resolution $(\CT,\pi_*,\pi^*)$ of a scheme $Y$ is {\sf strongly crepant}
if the relative Serre functor of $\CT$ over $\BD^b(\coh(Y))$ is isomorphic to the identity.
\end{definition}

Again, Grothendieck duality implies that for a geometric resolution crepancy 
and strong crepancy are equivalent. Moreover, it is not so difficult to show that
strong crepancy of a categorical resolution implies its weak crepancy, but the converse
is not true in general. To see this one can analyze the weak and strong crepancy
of categorical resolutions provided by Theorem~\ref{ldtocr}.

\begin{proposition}\label{crepcr}
In the setup of Theorem~$\ref{ldtocr}$ assume that $Y$ is Gorenstein
and $K_{\TY/Y} = (m-1)E$, where $m$ is the length of the left Lefschetz 
decomposition~\eqref{lde}. The corresponding categorical resolution 
$\TCC$ of $Y$ is weakly crepant if and only if 
\begin{equation}\label{wcrcond}
Li^*(Lf^*(\BD^\perf(Y))) \subset \CC_{m-1}.
\end{equation} 
Furthermore, $\TCC$ is strongly crepant if and only if~\eqref{lde} is rectangular, i.e.\
\begin{equation}\label{scrcond}
\CC_{m-1} = \dots = \CC_1 = \CC_0.
\end{equation} 
\end{proposition}

So, starting from a nonrectangular Lefschetz decomposition it is easy to produce
an example of a weakly crepant categorical resolution which is not strongly crepant.

\begin{example}
Take $X = Q^3 \subset \PP^4$ and let $Y$ be the cone over $X$ (i.e. a 4-dimensional quadratic cone).
Then the left Lefschetz collection 
\begin{equation*}
\BD^b(\coh(X)) = \langle \CC_2(-2), \CC_1(-1), \CC_0 \rangle
\qquad\text{with $\CC_0 = \langle \SSS, \CO_X \rangle$, $\CC_1 = \CC_2 = \langle \CO_X \rangle$}
\end{equation*}
($\SSS$ is the spinor bundle) gives a weakly crepant categorical resolution $\TCC$
of $Y$ which is not strongly crepant. In fact, if $q:\Tot_X(\CO_X(-1)) \to X$
is the canonical projection, the vector bundle $q^*\SSS$ is a spherical object
in $\TCC$ and the relative Serre functor is isomorphic to the corresponding spherical twist.
\end{example}

\subsection{Further questions}

Of course, the central question is Conjecture~\ref{mincr}. Theorem~\ref{ldtocr} shows
that it is closely related to the question of existence of minimal Lefschetz decompositions.

Another interesting question is to find new methods of construction
of minimal categorical resolutions. An interesting development in this direction
is the work \cite{Ab12} in which a notion of a {\sf wonderful resolution of singularities} 
(an analogue of wonderful compactifications) is introduced and it is shown
that a wonderful resolution gives rise to a weakly crepant
categorical resolution. This can be viewed as an advance on the first part 
of Proposition~\ref{crepcr}. It would be very interesting to find a generalization 
of the second part of this Proposition in the context of wonderful resolutions.

Another aspect is to find explicit constructions of minimal resolutions for interesting
varieties, such as Pfaffian varieties for example. Some of these arise naturally 
in the context of HP duality as we will see later.

\section{Examples of homological projective duality}\label{sec-ex}

If an HP duality for two varieties $X$ and $Y$ is proved, one gets as a consequence
an identification of the nontrivial components of the derived categories of linear 
sections of $X$ and $Y$.
Because of that it is clear that such a result is a very strong statement 
and is usually not so easy to prove. In this section we list several 
examples of HP duality. We assume that $\chr\kk = 0$ in this section.

\subsection{Linear duality}

Let $X = \PP_S(E)$ be a projectivization of a vector bundle $E$ on a scheme $S$
and assume that the map $f:X \to \PP(V)$ is {\em linear on fibers} of $X$ over~$S$.
In other words, we assume that $f$ is induced by an embedding of vector bundles
$E \to V\otimes\CO_S$ on $S$. In this case the line bundle $\CO_X(1) = f^*\CO_{\PP(V)}(1)$
is a Grothendieck line bundle for $X$ over $S$. By Theorem~\ref{thmorlsod}
we have a rectangular Lefschetz decomposition of $\BD^b(\coh(X))$ of length $m = \rk(E)$ with blocks
\begin{equation*}
\CA_0 = \CA_1 = \dots = \CA_{m-1} = p^*(\BD^b(\coh(S))),
\end{equation*}
where $p:X \to S$ is the projection. So, we are in the setup of HP duality and one can ask
what the dual variety is? 

The answer turns out to be given by a projectivization of 
another vector bundle over $S$.
Define $E^\perp$ as the kernel of the dual morphism 
\begin{equation*}
E^\perp := \Ker(V^\vee\otimes\CO_S \to E^\vee).
\end{equation*}
The projectivization $\PP_S(E^\perp)$ comes with a natural morphism
$g:\PP_S(E^\perp) \to \PP(V^\vee)$ and Theorem~\ref{thmorlsod} provides $\PP_S(E^\perp)$ with
a rectangular Lefschetz decomposition of length $N - m$ with blocks
$\CB_0 = \CB_1 = \dots = \CB_{N-m-1} = q^*(\BD^b(\coh(S)))$,
where $q:\PP_S(E^\perp) \to S$ is the projection.

\begin{theorem}[\cite{K07a}]
The projectivizations $X = \PP_S(E)$ and $Y = \PP_S(E^\perp)$ with their canonical
morphisms to $\PP(V)$ and $\PP(V^\vee)$ and the above Lefschetz decompositions
are homologically projectively dual to each other.
\end{theorem}

The picture visualizing this duality is very simple:
\begin{align*}
\ytableausetup{nosmalltableaux}
\begin{ytableau}
\none[\BD^b(\coh(S))]&\none&\none
&*(lightgray){\scriptstyle }
&*(lightgray)\scriptstyle 
&\none[\cdots]
&*(lightgray)\scriptstyle 
&\scriptstyle 
&\none[\cdots]
&\scriptstyle 
&\scriptstyle
&\scriptstyle 
&\scriptstyle 
&\none&\none&\none[\BD^b(\coh(S))]
\end{ytableau}
\end{align*}
with $m$ gray boxes and $N-m$ white boxes.

In the very special case of $S = \Spec\kk$ the bundle $E$ is just a vector space
and the variety $X$ is a (linearly embedded) projective subspace $\PP(E) \subset \PP(V)$.
Then the HP-dual variety is the orthogonal subspace $\PP(E^\perp) \subset \PP(V^\vee)$.
In particular, the dual of the space $\PP(V)$ itself with respect to its identity map
is the empty set.

\subsection{Quadrics}

There are two ways to construct a smooth quadric: one --- as a smooth hypersurface
of degree 2 in a projective space, and the other --- as a double covering
of a projective space ramified in a smooth quadric hypersurface.
These representations interchange in a funny way in HP duality.

Denote by $\SSS$ the spinor bundle on an odd dimensional quadric or 
one of the spinor bundles on the even dimensional quadric.

\begin{theorem}
$1)$ If $X = Q^{2m} \subset \PP^{2m+1}$ with Lefschetz decomposition given by
\begin{equation}\label{ldq2m}
\CA_0 = \CA_1 = \langle \SSS_X, \CO_X \rangle,
\qquad
\CA_2 = \CA_3 = \dots = \CA_{2m-1} = \langle \CO_X \rangle
\end{equation}
then the HP dual variety is the dual quadric $Y = Q^\vee \subset \check{\PP}^{2m+1}$
with the same Lefschetz decomposition.
\begin{align*}
\ytableausetup{smalltableaux}
\begin{ytableau}
\none[\CO_X]&\none
&*(lightgray)
&*(lightgray)
&*(lightgray)
&*(lightgray)
&\none&\none[\cdots]&\none
&*(lightgray)
&*(lightgray)
&&
&\none[\ \ \ \ \SSS_Y] \\
\none[\SSS_X]&\none
&*(lightgray)
&*(lightgray)
&&
&\none&\none[\cdots]&\none
&&&&&
\none[\ \ \ \ \CO_Y]&\none
\end{ytableau}
\end{align*}
$2)$ If $X = Q^{2m-1} \subset \PP^{2m}$ with Lefschetz decomposition given by
\begin{equation}\label{ldq2mm1}
\CA_0 = \langle \SSS_X, \CO_X \rangle,
\qquad
\CA_1 = \CA_2 = \dots = \CA_{2m-2} = \langle \CO_X \rangle
\end{equation}
then the HP dual variety is the double covering $Y \to \check{\PP}^{2m}$
ramified in the dual quadric $Q^\vee \subset \check{\PP}^{2m}$
with Lefschetz decomposition~\eqref{ldq2m}
\begin{align*}
\begin{ytableau}
\none[\CO_X]&\none
&*(lightgray)
&*(lightgray)
&*(lightgray)
&\none&\none[\cdots]&\none
&*(lightgray)
&*(lightgray)
&&
&\none[\ \ \ \ \SSS_Y] \\
\none[\SSS_X]&\none
&*(lightgray)
&&
&\none&\none[\cdots]&\none
&&&&&
\none[\ \ \ \ \CO_Y]&\none
\end{ytableau}
\end{align*}
$3)$ If $X = Q^{2m-1} \to \PP^{2m-1}$ is the double covering ramified 
in a quadric $\bar{Q} \subset \PP^{2m-1}$ with Lefschetz decomposition~\eqref{ldq2mm1}
then the HP dual variety is the double covering $Y \to \check{\PP}^{2m-1}$
ramified in the dual quadric $\bar{Q}^\vee \subset \check{\PP}^{2m-1}$
with the same Lefschetz decomposition.
\begin{align*}
\begin{ytableau}
\none[\CO_X]&\none
&*(lightgray)
&*(lightgray)
&*(lightgray)
&\none&\none[\cdots]&\none
&*(lightgray)
&*(lightgray)
&
&\none[\ \ \ \ \SSS_Y] \\
\none[\SSS_X]&\none
&*(lightgray)
&&
&\none&\none[\cdots]&\none
&&&&
\none[\ \ \ \ 	\CO_Y]&\none
\end{ytableau}
\end{align*}
\end{theorem}

\subsection{Veronese--Clifford duality}

Let $W$ be a vector space of dimension~$n$ and $V =S^2W$ its symmetric square.
We take $X = \PP(W)$ and consider its double Veronese embedding
$f:\PP(W) \to \PP(V)$. Then $f^*\CO_{\PP(V)}(1) = \CO_{\PP(W)}(2)$.
Beilinson's collection~\eqref{beiexcol} on $\PP(W)$ 
can be considered as a Lefschetz decomposition (with respect to $\CO_{\PP(W)}(2)$)
of $\BD^b(\coh(\PP(W)))$ with $\lfloor n/2 \rfloor$ blocks equal to
\begin{equation*}
\CA_0 = \CA_1 = \dots = \CA_{\lfloor n/2 \rfloor - 1} := \langle \CO_{\PP(W)}, \CO_{\PP(W)}(1) \rangle,
\end{equation*} 
and if $n$ is odd one more block
\begin{equation*}
\CA_{\lfloor n/2 \rfloor} := \langle \CO_{\PP(W)} \rangle.
\end{equation*} 
The universal hyperplane section $\CX$ of $X$ is nothing but the universal
quadric in $\PP(W)$ over the space $\PP(V^\vee) = \PP(S^2W^\vee)$ of all quadrics.
Then the quadratic fibration formula of Theorem~\ref{thm-qu-sod} gives an equivalence
of the HP dual category $\CC$ with the derived category $\BD^b(\coh(\PP(V^\vee),\Cliff_0))$
of coherent sheaves of modules over 
the even part of the universal Clifford algebra
\begin{equation*}
\Cliff_0 = \CO_{\PP(S^2W^\vee)} \oplus \Lambda^2W\otimes\CO_{\PP(S^2W^\vee)}(-1) \oplus \dots,
\end{equation*}
on the space $\PP(S^2W^\vee)$ of quadrics. We will consider the pair $(\PP(S^2W^\vee),\Cliff_0)$
as a noncommutative variety and call it {\sf the Clifford space}.

\begin{theorem}[Veronese--Clifford duality, \cite{K08a}]
The HP dual of the projective space $X = \PP(W)$ in the double Veronese embedding $\PP(W) \to \PP(S^2W)$
is the Clifford space $Y = (\PP(S^2W^\vee),\Cliff_0)$. 
\end{theorem}

The HP dual Lefschetz decomposition of the Clifford space is given by the full exceptional collection
\begin{equation*}
\BD^b(\coh(\PP(S^2W^\vee,\Cliff_0))) = 
\langle \Cliff_{1-n^2},\Cliff_{2-n^2},\dots,\Cliff_{-1},\Cliff_{0} \rangle,
\end{equation*}
where 
\begin{equation*}
\Cliff_1 = W\otimes\CO_{\PP(S^2W^\vee)} \oplus \Lambda^3W\otimes\CO_{\PP(S^2W^\vee)}(-1) \oplus \dots
\end{equation*}
is the odd part of the Clifford algebra and $\Cliff_{k-2} = \Cliff_k \otimes \CO_{\PP(S^2W^\vee)}(-1)$
for each $k \in \ZZ$. The picture visualizing this duality is:
\begin{align*}
\ytableausetup{smalltableaux}
\begin{array}{c}
\begin{ytableau}
\none[\CO]&\none
&*(lightgray)
&*(lightgray)
&\none&\none[\cdots]&\none
&*(lightgray)
&\scriptstyle 
&\none&\none[\cdots]&\none
&\scriptstyle
&\scriptstyle 
&\scriptstyle 
&\none&\none[\Cliff_1]\\
\none[\CO(1)]&\none
&*(lightgray)
&*(lightgray)
&\none&\none[\cdots]&\none
&*(lightgray)
&\scriptstyle 
&\none&\none[\cdots]&\none
&\scriptstyle
&\scriptstyle
&\scriptstyle 
&\none&\none[\Cliff_0]
\end{ytableau}\\[2ex]
\text{$n$ even}
\end{array}
&\text{\quad and\ }&
\begin{array}{c}
\begin{ytableau}
\none[\CO]&\none&
*(lightgray)
&*(lightgray)
&\none&\none[\cdots]&\none
&*(lightgray)
&*(lightgray)
&\scriptstyle 
&\none&\none[\cdots]&\none
&\scriptstyle 
&\scriptstyle
&\scriptstyle 
&\none&\none[\Cliff_1]\\
\none[\CO(1)]&\none&
*(lightgray)
&*(lightgray)
&\none&\none[\cdots]&\none
&*(lightgray)
&\scriptstyle 
&\scriptstyle
&\none&\none[\cdots]&\none
&\scriptstyle 
&\scriptstyle
&\scriptstyle 
&\none&\none[\Cliff_0]
\end{ytableau}\\[2ex]
\text{$n$ odd}
\end{array}
\end{align*}
for even $n$ it has $n/2$ gray columns and $n^2/2$ white columns, 
and for odd $n$ it has $(n-1)/2$ gray columns, one mixed column, and $(n^2-1)/2$ white columns.

\subsection{Grassmannian--Pfaffian duality}\label{secgpd}

The most interesting series of examples is provided by Grassmannians $\Gr(2,m)$
of two-dimensional subspaces in an $m$-dimensional vector space.

Let $W$ be a vector space of dimension $m$ and let $V = \Lambda^2W$, the space of bivectors.
The group $\GL(W)$ acts on the projective space $\PP(\Lambda^2W)$ with orbits
indexed by the rank of a bivector which is always even and ranges 
from $2$ to $2\lfloor m/2\rfloor$. We denote by $\Pf(2k,W)$ the closure of the orbit 
consisting of bivectors of rank $2k$ and call it the {\sf $k$-th Pfaffian variety}.
Clearly, the smallest orbit $\Pf(2,W)$ is smooth and coincides with the Grassmannian
$\Gr(2,W)$ in its Pl\"ucker embedding. Another smooth Pfaffian variety is the maximal one ---
$\Pf(2\lfloor m/2 \rfloor,W) = \PP(\Lambda^2W)$. All the intermediate Pfaffians
are singular with
$\sing(\Pf(2k,W)) = \Pf(2k-2,W)$.
The submaximal Pfaffian variety $\Pf(2\lfloor m/2\rfloor - 2,W^\vee)$ of the dual space is
classically projectively dual to the Grassmannian $\Gr(2,W)$. 
This suggests a possible HP duality between them. 

To make a precise statement we should choose a Lefschetz decomposition of $\BD^b(\coh(\Gr(2,W)))$.
A naive choice is to take Kapranov's collection~\eqref{kapexcol}.
It can be considered as a Lefschetz decomposition on $X := \Gr(2,m)$ with $m-1$ blocks
\begin{equation*}
\CA_0 = \langle \CO_X,\CU_X^\vee,\dots,S^{m-2}\CU_X^\vee \rangle,\  
\CA_1 = \langle \CO_X,\CU_X^\vee,\dots,S^{m-3}\CU_X^\vee \rangle,\ \dots,\ 
\CA_{m-2} = \langle \CO_X \rangle.
\end{equation*}
However, it is very far from being minimal. It turns out that a reasonable result
can be obtained for another Lefschetz decomposition
\begin{equation*}
\BD^b(\coh(\Gr(2,m)) = \langle \CA_0,\CA_1(1),\dots,\CA_{m-1}(m-1) \rangle
\end{equation*} 
with 
\begin{equation}\label{ldgodd}
\CA_0 = \dots = \CA_{m-1} = \langle \CO_X,\CU_X^\vee,\dots,S^{(m-1)/2}\CU_X^\vee \rangle.
\end{equation} 
if $m$ is odd, and with
\begin{align}
\CA_0 = \dots = \CA_{m/2-1}  = &\ \langle \CO_X,\CU_X^\vee,\dots,S^{m/2-1}\CU_X^\vee \rangle,\nonumber\\
\CA_{m/2} = \dots = \CA_{m-1} = &\ \langle \CO_X,\CU_X^\vee,\dots,S^{m/2-2}\CU_X^\vee \rangle,\label{ldgeven}
\end{align} 
if $m$ is even.

\begin{conjecture}\label{conjgrpf}
The HP dual of the Grassmannian $\Gr(2,W)$ with Lefschetz decomposition~\eqref{ldgeven} 
{\rm(}or~\eqref{ldgodd} depending on the parity of $m = \dim W${\rm)} is given 
by a {minimal} categorical resolution of the submaximal Pfaffian $\Pf(2\lfloor m/2 \rfloor - 2,W^\vee)$. 
When $m$ is odd, this resolution is strongly crepant.
\end{conjecture}

This conjecture is proved for $m \le 7$ in~\cite{K06b}. 
In fact, for $m = 2$ and $m = 3$ one has $\Gr(2,W) = \PP(\Lambda^2W)$
and linear duality applies. For $m = 4$ and $m = 5$ the submaximal Pfaffian $\Pf(2,W^\vee)$
coincides with the Grassmannian, and the above duality is the duality for Grassmannians:
\begin{align*}
\ytableausetup{smalltableaux}
\begin{array}{c}
\Gr(2,4)\\ 
\begin{ytableau}
\none[\scriptscriptstyle\CO_X]&\none
&*(lightgray)&*(lightgray)&*(lightgray)&*(lightgray)
&&&
\none&\none[\scriptscriptstyle \CU_Y]
\\
\none[\scriptscriptstyle\CU^\vee_X]&\none
&*(lightgray)&*(lightgray)
&&&&&
\none&\none[\scriptscriptstyle \CO_Y]
\end{ytableau}
\end{array}
\qquad
\qquad
\begin{array}{c}
\Gr(2,5)\\
\begin{ytableau}
\none[\scriptscriptstyle\CO_X]&\none
&*(lightgray)&*(lightgray)&*(lightgray)&*(lightgray)&*(lightgray)
&&&&&&
\none&\none[\scriptscriptstyle \CU_Y]
\\
\none[\scriptscriptstyle\CU^\vee_X]&\none
&*(lightgray)&*(lightgray)&*(lightgray)&*(lightgray)&*(lightgray)
&&&&&&
\none&\none[\scriptscriptstyle \CO_Y]
\end{ytableau}
\end{array}
\end{align*}
For $m = 6$ and $m = 7$ the submaximal Pfaffian $Y = \Pf(4,W^\vee)$ is singular, but its appropriate
categorical resolutions can be constructed by Theorem~\ref{ldtocr}. 
It turns out that these
resolutions indeed are HP dual to the corresponding Grassmannians:
\begin{align*}
\ytableausetup{smalltableaux}
\begin{array}{c}
\Gr(2,6)\qquad\qquad\Pf(4,6)\\
\begin{ytableau}
*(lightgray)&*(lightgray)&*(lightgray)&*(lightgray)&*(lightgray)&*(lightgray)
&&&&&&&&&
\\
*(lightgray)&*(lightgray)&*(lightgray)&*(lightgray)&*(lightgray)&*(lightgray)
&&&&&&&&&
\\
*(lightgray)&*(lightgray)&*(lightgray)
&&&&&&&&&&&&
\end{ytableau}
\end{array}
\qquad
\begin{array}{c}
\Gr(2,7)\qquad\qquad\qquad\qquad\Pf(4,7)\\
\begin{ytableau}
*(lightgray)&*(lightgray)&*(lightgray)&*(lightgray)&*(lightgray)&*(lightgray)&*(lightgray)
&&&&&&&&&&&&&&
\\
*(lightgray)&*(lightgray)&*(lightgray)&*(lightgray)&*(lightgray)&*(lightgray)&*(lightgray)
&&&&&&&&&&&&&&
\\
*(lightgray)&*(lightgray)&*(lightgray)&*(lightgray)&*(lightgray)&*(lightgray)&*(lightgray)
&&&&&&&&&&&&&&
\end{ytableau}
\end{array}
\end{align*}
For $m \ge 8$ this construction of a categorical resolution does not work.
However it is plausible that the Pfaffians have wonderful resolutions
of singularities, so a development of~\cite{Ab12} may solve the question.

\subsection{The spinor duality}

Let $W$ be a vector space of even dimension $2m$ and $q \in S^2W^\vee$ a nondegenerate quadratic form.
The isotropic Grassmannian of $m$-dimensional subspaces in $W$ has two connected components,
abstractly isomorphic to each other and called {\sf spinor varieties} $\SSSS_m$. These are homogeneous
spaces of the spin group $\Spin(W)$ with the embedding into $\PP(\Lambda^mW)$
given by the square of the generator of the Picard group, while the generator itself gives an embedding
into the projectivization $\PP(V)$ of a half-spinor representation $V$ (of dimension $2^{m-1}$) of $\Spin(W)$. 
For small $m$ the spinor varieties
are very simple (because the spin-group simplifies): in fact, $\SSSS_1$ is a point,
$\SSSS_2 = \PP^1$, $\SSSS_3 = \PP^3$, and $\SSSS_4 = Q^6$. The first interesting
example is $\SSSS_5$. 

\begin{theorem}[\cite{K06a}]\label{hpdspinor}
The spinor variety $X = \SSSS_5$ has a Lefschetz decomposition 
\begin{equation}\label{ldspin5}
\BD^b(\coh(X)) = \langle \CA_0,\CA_1(1),\dots,\CA_7(7) \rangle
\quad\text{with}\quad
\CA_0 = \dots = \CA_7 = \langle \CO_X, \CU_5^\vee \rangle.
\end{equation}
The HP dual variety is the same spinor variety $Y = \SSSS_5$.
\begin{align*}
\ytableausetup{smalltableaux}
\begin{array}{c}
\SSSS_5\qquad\qquad\qquad\qquad\SSSS_5\\
\begin{ytableau}
*(lightgray)&*(lightgray)&*(lightgray)&*(lightgray)&*(lightgray)&*(lightgray)&*(lightgray)&*(lightgray)&
&&&&&&&
\\
*(lightgray)&*(lightgray)&*(lightgray)&*(lightgray)&*(lightgray)&*(lightgray)&*(lightgray)&*(lightgray)&
&&&&&&&
\end{ytableau}
\end{array}
\end{align*}
\end{theorem}

\subsection{Incomplete dualities}

It is often quite hard to give a full description of the HP dual variety.
On the other hand, there is sometimes an open dense subset $U \subset \PP(V^\vee)$ for which 
there is a description of the category $\CC_U$ obtained by a base change $U \to \PP(V^\vee)$
from the HP dual category $\CC$. If $Y_U$ is a (noncommutative) variety such that $Y_U \cong \CC_U$
(a $U$-linear equivalence), we will say that $Y_U$ is the {\sf HP dual of $X$ over $U$},
or an {\sf incomplete HP dual variety}.

\begin{proposition}
If $Y_U$ is an HP dual of $X$ over an open subset $U \subset \PP(V^\vee)$ then the semiorthogonal
decompositions~\eqref{dxl} and~\eqref{dyl} hold for any subspace $L \subset V^\vee$ such that
$\PP(L) \subset U$ and the varieties $X_L$ and $Y_L$
have expected dimensions.
\end{proposition}

There are several examples of HP duality when only an incomplete dual variety is known.
Below we discuss two of them.
Let $W$ be a vector space of dimension 6 with a symplectic form $\omega$,
and $X = \LGr(3,W) \subset \PP^{13}$ the corresponding Lagrangian Grassmannian. 
The classical projectively dual variety $Y_1 := X^\vee \subset \check{\PP}^{13}$ is a quartic
hypersurface which is singular along a $9$-dimensional variety $Y_2 \subset Y_1$.
One can check that $X$ has a rectangular Lefschetz decomposition (see~\cite{S01})
\begin{align}
&\BD^b(\coh(X)) = \langle \CA_0,\CA_1(1),\CA_2(2),\CA_3(3) \rangle
\quad\text{with}\quad\label{ldlg6}\\
&\CA_0 = \CA_1 = \CA_2 = \CA_3 = \langle \CO_X, \CU_X^\vee \rangle.\nonumber
\end{align} 
The HP dual variety in this case is described only
over the open set $U = \check{\PP}^{13} \setminus Y_2$, see~\cite{K06a}.
For this a morphism $\pi:\tilde{Y} \to Y_1$ which is a nondegenerate
conic bundle over $Y_1\setminus Y_2$ is constructed in {\em loc.\ cit.}
Let $\CR$ be the associated quaternion algebra.

\begin{theorem}[\cite{K06a}]\label{hpdlg36}
If $X = \LGr(3,6) \subset \PP^{13}$ with Lefschetz decomposition~\eqref{ldlg6}
then the noncommutative variety $(Y_1\setminus Y_2,\CR)$ is the HP dual of $X$
over $\check{\PP}^{13} \setminus Y_2$.
\begin{align*}
\begin{array}{c}
\LGr(3,6)\qquad\quad (Y_1\setminus Y_2,\CR)\\
\ytableausetup{smalltableaux}
\begin{ytableau}
*(lightgray)&*(lightgray)&*(lightgray)&*(lightgray)&
&&&&&&&&&
\\
*(lightgray)&*(lightgray)&*(lightgray)&*(lightgray)&
&&&&&&&&&
\end{ytableau}
\end{array}
\end{align*}
\end{theorem}

\bigskip

Another example is related to the simple algebraic group $\Gt$.
Let $X$ be the orbit of the highest
vector in the projectivization $\PP^{13} = \PP(V)$ of the adjoint representation $V$ of $\Gt$. 
Then $X$ can also be realized as the zero locus of a global section of the vector bundle 
$\Lambda^4(W/\CU)$ on the Grassmannian $\Gr(2,W)$ for a 7-dimensional fundamental 
representation $W$ of $\Gt$ corresponding to a generic 3-form $\lambda \in \Lambda^3W^\vee$. 
By that reason we use the notation $\Gtgr(2,W)$ for $X$.

The classical projectively dual variety $Y_1 := X^\vee \subset \check{\PP}^{13}$ is a sextic
hypersurface which is singular along a $10$-dimensional variety $Y_2 \subset Y_1$.
One can check that $X$ has a rectangular Lefschetz decomposition 
\begin{equation}\label{ldgtgr}
\BD^b(\coh(X)) = \langle \CA_0,\CA_1(1),\CA_2(2) \rangle
\quad\text{with}\quad
\CA_0 = \CA_1 = \CA_2 = \langle \CO_X, \CU_X^\vee \rangle.
\end{equation} 
As in the previous case, the HP dual variety is described only over the open set
$U = \check{\PP}^{13} \setminus Y_2$, see~\cite{K06a}.
A morphism $\pi:\tilde{Y} \to Y_0$ to the double covering $Y_0 \to \check{\PP}^{13}$ 
ramified in $Y_1$ which is a Severi--Brauer variety with fiber $\PP^2$ over $\TY \setminus Y_2$
is constructed in {\em loc.\ cit.}
Let $\CR$ be the associated Azumaya algebra on $Y_0 \setminus Y_2$.

\begin{theorem}[\cite{K06a}]\label{hpdg2}
If $X = \Gtgr(2,7)$ with Lefschetz decomposition~\eqref{ldgtgr}
then the noncommutative variety $(Y_0 \setminus Y_2,\CR)$ is the HP dual of $X$
over $\check{\PP}^{13} \setminus Y_2$.
\begin{align*}
\begin{array}{c}
\Gtgr(2,7)\qquad (Y_0 \setminus Y_2,\CR)\\
\ytableausetup{smalltableaux}
\begin{ytableau}
*(lightgray)&*(lightgray)&*(lightgray)&
&&&&&&&&&&
\\
*(lightgray)&*(lightgray)&*(lightgray)&
&&&&&&&&&&
\end{ytableau}
\end{array}
\end{align*}
\end{theorem}

\subsection{Conjectures}

Note that $\Gr(2,W) = \Pf(2,W)$. This suggests a generalization of the Grassmannian--Pfaffian duality 
to higher Pfaffians.

\begin{conjecture}\label{conjpfpf}
For any $k$ there is an HP duality between appropriate {minimal} categorical 
resolutions of the Pfaffians $\Pf(2k,W)$ and $\Pf(2(\lfloor m/2\rfloor-k),W^\vee)$. 
When $m = \dim W$ is odd these resolutions are strongly crepant.
\end{conjecture}

Below are the expected pictures for the HP duality for $\Pf(4,8)$ and for $\Pf(4,9)$:
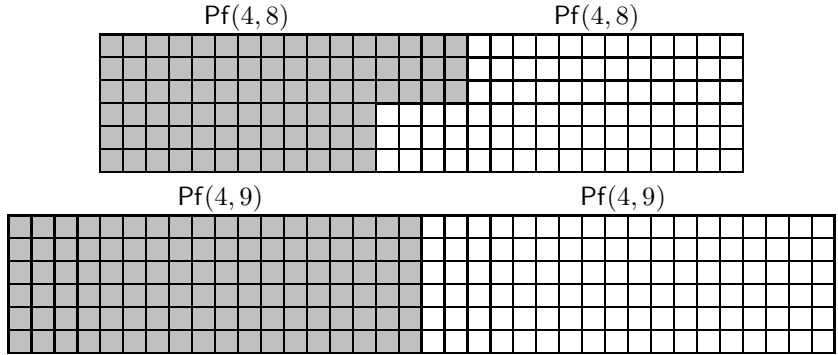
\begin{align*}
\ytableausetup{smalltableaux}
\begin{array}{c}
\begin{array}{c}
\Pf(4,8)
\qquad\qquad\qquad\qquad\qquad
\Pf(4,8)\\
\begin{ytableau}
*(lightgray)&*(lightgray)&*(lightgray)&*(lightgray)&*(lightgray)&*(lightgray)&*(lightgray)&*(lightgray)
&*(lightgray)&*(lightgray)&*(lightgray)&*(lightgray)&*(lightgray)&*(lightgray)&*(lightgray)&*(lightgray)
&&&&&&&&&&&&
\\
*(lightgray)&*(lightgray)&*(lightgray)&*(lightgray)&*(lightgray)&*(lightgray)&*(lightgray)&*(lightgray)
&*(lightgray)&*(lightgray)&*(lightgray)&*(lightgray)&*(lightgray)&*(lightgray)&*(lightgray)&*(lightgray)
&&&&&&&&&&&&
\\
*(lightgray)&*(lightgray)&*(lightgray)&*(lightgray)&*(lightgray)&*(lightgray)&*(lightgray)&*(lightgray)
&*(lightgray)&*(lightgray)&*(lightgray)&*(lightgray)&*(lightgray)&*(lightgray)&*(lightgray)&*(lightgray)
&&&&&&&&&&&&
\\
*(lightgray)&*(lightgray)&*(lightgray)&*(lightgray)&*(lightgray)&*(lightgray)
&*(lightgray)&*(lightgray)&*(lightgray)&*(lightgray)&*(lightgray)&*(lightgray)
&&&&&&&&&&&&&&&&
\\
*(lightgray)&*(lightgray)&*(lightgray)&*(lightgray)&*(lightgray)&*(lightgray)
&*(lightgray)&*(lightgray)&*(lightgray)&*(lightgray)&*(lightgray)&*(lightgray)
&&&&&&&&&&&&&&&&
\\
*(lightgray)&*(lightgray)&*(lightgray)&*(lightgray)&*(lightgray)&*(lightgray)
&*(lightgray)&*(lightgray)&*(lightgray)&*(lightgray)&*(lightgray)&*(lightgray)
&&&&&&&&&&&&&&&&
\end{ytableau}
\end{array}\\[7ex]
\ytableausetup{smalltableaux}
\begin{array}{c}
\Pf(4,9)
\qquad\qquad\qquad\qquad\qquad\qquad
\Pf(4,9)\\
\begin{ytableau}
*(lightgray)&*(lightgray)&*(lightgray)&*(lightgray)&*(lightgray)&*(lightgray)&*(lightgray)&*(lightgray)&*(lightgray)
&*(lightgray)&*(lightgray)&*(lightgray)&*(lightgray)&*(lightgray)&*(lightgray)&*(lightgray)&*(lightgray)&*(lightgray)
&&&&&&&&&&&&&&&&&&
\\
*(lightgray)&*(lightgray)&*(lightgray)&*(lightgray)&*(lightgray)&*(lightgray)&*(lightgray)&*(lightgray)&*(lightgray)
&*(lightgray)&*(lightgray)&*(lightgray)&*(lightgray)&*(lightgray)&*(lightgray)&*(lightgray)&*(lightgray)&*(lightgray)
&&&&&&&&&&&&&&&&&&
\\
*(lightgray)&*(lightgray)&*(lightgray)&*(lightgray)&*(lightgray)&*(lightgray)&*(lightgray)&*(lightgray)&*(lightgray)
&*(lightgray)&*(lightgray)&*(lightgray)&*(lightgray)&*(lightgray)&*(lightgray)&*(lightgray)&*(lightgray)&*(lightgray)
&&&&&&&&&&&&&&&&&&
\\
*(lightgray)&*(lightgray)&*(lightgray)&*(lightgray)&*(lightgray)&*(lightgray)&*(lightgray)&*(lightgray)&*(lightgray)
&*(lightgray)&*(lightgray)&*(lightgray)&*(lightgray)&*(lightgray)&*(lightgray)&*(lightgray)&*(lightgray)&*(lightgray)
&&&&&&&&&&&&&&&&&&
\\
*(lightgray)&*(lightgray)&*(lightgray)&*(lightgray)&*(lightgray)&*(lightgray)&*(lightgray)&*(lightgray)&*(lightgray)
&*(lightgray)&*(lightgray)&*(lightgray)&*(lightgray)&*(lightgray)&*(lightgray)&*(lightgray)&*(lightgray)&*(lightgray)
&&&&&&&&&&&&&&&&&&
\\
*(lightgray)&*(lightgray)&*(lightgray)&*(lightgray)&*(lightgray)&*(lightgray)&*(lightgray)&*(lightgray)&*(lightgray)
&*(lightgray)&*(lightgray)&*(lightgray)&*(lightgray)&*(lightgray)&*(lightgray)&*(lightgray)&*(lightgray)&*(lightgray)
&&&&&&&&&&&&&&&&&&
\end{ytableau}
\end{array}
\end{array}
\end{align*}

\bigskip

It is also expected that there is a generalization of the Veronese--Clifford duality.
Consider a vector space $W$ of dimension $n$ and its symmetric square $V = S^2W$.
The group $\GL(W)$ acts on the projective space $\PP(V) = \PP(S^2W)$ with orbits indexed
by the rank of a tensor. 
We denote by $\Sigma(k,W) \subset \PP(S^2W)$ 
the closure of the orbit consisting of symmetric tensors of rank~$k$. The smallest orbit $\Sigma(1,W)$
is smooth and coincides with the double Veronese embedding of the projective space $\PP(W)$.
On the other hand, $\Sigma(n,W) = \PP(S^2W)$ is also smooth. All the intermediate varieties $\Sigma(k,W)$
are singular with
$\sing(\Sigma(k,W)) = \Sigma(k-1,W)$.
The classical projective duality acts on these varieties 
by $\Sigma(k,W)^\vee = \Sigma(n-k,W^\vee)$. However, the HP duality
is organized in a much more complicated way. 

Besides $\Sigma(k,W)$ itself one can consider its modifications: 
\begin{itemize}
\item 
the Clifford modification $(\Sigma(k,W),\Cliff_0)$ for the natural sheaf 
of even parts of Clifford algebras on it, and
\item 
(for even $k$) the double covering $\tilde\Sigma(k,W)$ of $\Sigma(k,W)$
corresponding to the central subalgebra in $\Cliff_0$ as in~\eqref{stilde},
\end{itemize}
and their minimal categorical resolutions.
It seems that HP duality interchanges in a complicated way modifications of different type.
For example, besides the original Veronese--Clifford duality between $\Sigma(1,n)$
and $(\Sigma(n,n),\Cliff_0)$ there are strong indications that
(the minimal resolution of) $\Sigma(2,4)$ is HP dual to
(the minimal resolution of) the double covering $\tilde\Sigma(4,4)$ 
of $\Sigma(4,4) = \PP^9$ (see~\cite{IK10}),
(the minimal resolution of) $\Sigma(2,5)$ is HP dual to
(the minimal resolution of) the double covering $\tilde\Sigma(4,5)$ 
of $\Sigma(4,5)$ (see~\cite{HT}), while
(the minimal resolution of) the double covering $\tilde\Sigma(2,n)$ of $\Sigma(2,n)$ is HP dual to
(the minimal resolution of) $\Sigma(n-1,n)$ for all $n$ (this can be easily deduced
from the linear duality).

\section{Varieties of small dimension}\label{sec-var}

Let us list what is known about semiorthogonal decompositions of smooth projective 
varieties by dimension. In this section we assume that $\kk = \mathbb{C}$.

\subsection{Curves}

Curves are known to have no nontrivial semiorthogonal decompositions 
with the only exception of $\PP^1$ (for which every semiorthogonal decomposition
coincides with the Beilinson decomposition up to a twist), see~\cite{Ok11}.

\subsection{Surfaces}\label{ssec-surfaces}

For surfaces the situation is more complicated. Of course, by the blowup formula
any surface has a semiorthogonal decomposition with several exceptional objects
and the derived category of a minimal surface as components. In particular,
any rational surface has a full exceptional collection. Moreover, for $\PP^2$ 
it is known that any full exceptional collection can be obtained from Beilinson's 
collection by mutations (which are related to Markov numbers and toric degenerations 
of $\PP^2$). For other del Pezzo surfaces all exceptional objects have been classified \cite{KO}, 
and moreover, three-blocks exceptional collections were constructed \cite{KN}, but the complete 
picture is not known.

For minimal ruled surfaces, of course there is a semiorthogonal decomposition 
into two copies of the derived category of the curve which is the base of the ruling.

For surfaces of Kodaira dimension 0 it is well known that there are no nontrivial semiorthogonal
decompositions for K3 and abelian surfaces. For Enriques surfaces there may be an exceptional
collection of line bundles of length up to 10 (see~\cite{Z}), and for so-called nodal Enriques surfaces 
the complementary component is related to the Artin--Mumford quartic double solid \cite{IK10}.
For bielliptic surfaces nothing is known.

For surfaces of Kodaira dimension 1 with globally generated canonical class (and even for those
without multiple fibers and $p_g > 0$) there are no semiorthogonal decompositions by~\cite{KaOk12}.
For others nothing is known as well.


Finally, for surfaces of general type there is an unexpectedly rich theory of semiorthogonal
decompositions. In fact, for many surfaces of general type with $p_g = q = 0$ (the classical 
Godeaux surface, the Beauville surface, the Burniat surfaces, the determinantal Barlow surface, 
some fake projective planes)
exceptional collections of length equal to the rank of the Grothendieck group have been constructed
in \cite{BBS13,GS13,AO13,BBKS12,GKMS13,F13}.
The collections, however, are not full. The complementary components have
finite (or even zero) Grothendieck group and trivial Hochschild homology and by that reason
they are called {\sf quasiphantom} or {\sf phantom} categories. The phantoms cannot be detected 
by additive invariants, but one can use Hochschild cohomology instead, see~\cite{K12}.

An interesting feature here is that the structure of the constructed exceptional collections resembles 
very much the structure of exceptional collections of del Pezzo surfaces with the same $K^2$. The only
(but a very important) difference is that whenever there is a $\Hom$-space between exceptional
bundles on del Pezzo, the corresponding exceptional bundles on the surface of general type
have $\Ext^2$-space. This seemingly small difference, however, has a very strong effect
on the properties of the category. See more details in {\em loc.\ cit.}

\subsection{Fano 3-folds}

For derived categories of threefolds (and higher dimensional varieties) there are no classification results
(as there is no classification of threefolds). Of course, as it already was mentioned for varieties 
with trivial (or globally generated) canonical class there are no nontrivial decompositions. So, from now on 
we will discuss Fano varieties. 

In dimension 3 all Fano varieties were classified in the works of Fano, Iskovskikh and Mukai.
All Fano 3-folds with Picard number greater than 1 are either the blowups of other Fano varieties
with centers in points and smooth curves (and then their derived category reduces
to the derived category of a Fano 3-fold with smaller Picard number), 
or conic bundles over rational surfaces (see Tables 12.3--12.6 of~\cite{IP}). 
For conic bundles one can use the quadratic bundle formula (Theorem~\ref{thm-qu-sod}). 
It gives a semiorthogonal decomposition with several exceptional objects and the derived 
category of sheaves of modules over the even part of the Clifford algebra on the base of the bundle.

If the Picard number is 1, the next discrete invariant of a Fano 3-fold to look at is the {\sf index},
i.e.\ the maximal integer dividing the canonical class. By Fujita's Theorem 
the only Fano 3-folds of index greater than 2 are $\PP^3$ and $Q^3$.
Their derived categories are well understood,
so let us turn to 3-folds of index 2 and 1.

For a Fano 3-fold $Y$ of index 2 the pair of line bundles $(\CO_Y,\CO_Y(1))$ is exceptional
and gives rise to a semiorthogonal decomposition
\begin{equation}\label{sodi2}
\BD^b(\coh(Y)) = \langle \CB_Y, \CO_Y, \CO_Y(1) \rangle.
\end{equation}
The component $\CB_Y$ is called the {\sf nontrivial component} of $\BD^b(\coh(Y))$. 

A similar decomposition can be found for a Fano 3-fold $X$ of index 1 and degree $d_X := (-K_X)^3$
which is not divisible by 4 (the degree of a 3-fold of index 1 is always even).
By a result of Mukai~\cite{Mu92} on such $X$ there is an exceptional 
vector bundle $\CE_X$ of rank $2$ with $c_1(\CE_X) = K_X$, which is moreover orthogonal 
to the structure sheaf of~$X$. In other words, $(\CE_X,\CO_X)$ is an exceptional pair
and there is a semiorthogonal decomposition
\begin{equation}\label{sodi1}
\BD^b(\coh(X)) = \langle \CA_X, \CE_X, \CO_X \rangle.
\end{equation}
The component $\CA_X$ is called the {\sf nontrivial component} of $\BD^b(\coh(X))$. 

It is rather unexpected that the nontrivial parts $\CB_Y$ and $\CA_X$
for a Fano 3-fold $Y$ of index $2$ and degree $d_Y := (-K_Y/2)^3$ 
and for a Fano 3-fold $X$ of index 1 and degree $d_X = 4d_Y + 2$
have the same numerical characteristics,
and are, moreover, expected to belong to the same deformation family of categories.
In fact, this expectation is supported by the following result. Recall that the degree
of a Fano 3-fold of index 2 with Picard number 1 satisfies $1 \le d \le 5$, while the degree
of a Fano 3-fold of index 1 with Picard number 1 is even and satisfies $2 \le d \le 22$, $d \ne 20$.
So there are actually 5 cases to consider.

\begin{theorem}[\cite{K09a}]
For $3 \le d \le 5$ each category $\CB_{Y_d}$ is equivalent to some category $\CA_{X_{4d+2}}$
and vice versa.
\end{theorem}

See {\em loc.\ cit.}\/ for a precise statement.
In fact, for $d = 5$ the category is rigid and is equivalent to the derived category
of representations of the quiver with 2 vertices and 3 arrows from the first vertex to the second
(this follows from the construction of explicit exceptional collections in the derived 
categories of $Y_5$ and $X_{22}$, see~\cite{O91} and~\cite{K96}). 
Further, for $d = 4$ each of the categories $\CB_{Y_4}$
and $\CA_{X_{18}}$ is equivalent to the derived category of a curve of genus 2, and moreover,
each smooth curve appears in both pictures. This follows from HP duality for the double Veronese embedding
of $\PP^5$ and from HP duality for $\Gt$ Grassmannian respectively \cite{K06a}. Finally, for $d = 3$
no independent description of the category in question is known, but the HP duality 
for the Grassmannian $\Gr(2,6)$ gives the desired equivalence (see~\cite{K04,K06b}).

It turns out, however, that already for $d = 2$ the situation is more subtle. 
It seems that in that case the categories $\CB_{Y_2}$ lie at the boundary of the family
of categories $\CA_{X_{10}}$. And for $d = 1$ the situation is completely unclear.

The situation with Fano 3-folds of index 1 and degree divisible by~4 is somewhat different.
For such threefolds it is, in general, not clear how one can construct an exceptional pair.
However, for $d_X = 12$ and $d_X = 16$ this is possible. For $d_X = 12$ Mukai has proved~\cite{Mu92}
that there is an exceptional pair $(\CE_5,\CO_X)$ where $\CE_5$ is a rank 5 exceptional bundle
with $c_1(\CE_5) = 2K_X$. Using HP duality for the spinor variety~$\SSSS_5$ one can check 
that this pair extends to a semiorthogonal decomposition
\begin{equation*}
\BD^b(\coh(X_{12})) = \langle \BD^b(\coh(C_7)), \CE_5, \CO_{X_{12}} \rangle,
\end{equation*}
where $C_7$ is a smooth curve of genus $7$ (see~\cite{K05,K06a}). 
Analogously, for $d_X = 16$ Mukai has constructed~\cite{Mu92} an exceptional
bundle $\CE_3$ of rank $3$ with $c_1(\CE_3) = K_X$. Using HP duality for $\LGr(3,6)$ one can check
that there is a semiorthogonal decomposition
\begin{equation*}
\BD^b(\coh(X_{16})) = \langle \BD^b(\coh(C_3)), \CE_3, \CO_{X_{16}} \rangle,
\end{equation*}
where $C_3$ is a smooth curve of genus $3$ \cite{K06a}.

\subsection{Fourfolds}

Of course, for Fano 4-folds we know much less than for 3-folds. So, we will not even
try to pursue a classification, but will restrict attention to some very special cases of interest.

Maybe one of the most interesting 4-folds is the cubic 4-fold. One of its salient
features is the hyperk\"ahler structure on the Fano scheme of lines, which turns out
to be a deformation of the second Hilbert scheme of a K3 surface. This phenomenon
has a nice explanation from the derived categories point of view.

\begin{theorem}[\cite{K10}]
Let $Y \subset \PP^5$ be a cubic $4$-fold. Then there is a semiorthogonal decomposition
\begin{equation*}
\BD^b(\coh(Y)) = \langle \CA_Y,\CO_Y,\CO_Y(1),\CO_Y(2) \rangle,
\end{equation*}
and its nontrivial component $\CA_Y$ is a Calabi--Yau category of dimension $2$.
Moreover, $\CA_Y$ is equivalent to the derived category of coherent sheaves 
on a K$3$ surface, at least if
$Y$ is a Pfaffian cubic $4$-fold, or
if $Y$ contains a plane $\Pi$ and a $2$-cycle $Z$ such that 
$\deg Z + Z \cdot \Pi \equiv 1 \bmod 2$.
\end{theorem}

To establish this result for Pfaffian cubics one can use HP duality for $\Gr(2,6)$.
The associated K3 is then a linear section of this Grassmannian.
For cubics with a plane a quadratic bundle formula for the projection of $Y$ from the plane $\Pi$
gives the result.
The K3 surface then is the double covering of $\PP^2$ ramified in a sextic curve, 
and the cycle $Z$ gives a splitting of the requisite Azumaya algebra on this~K3.

For generic $Y$ the category $\CA_Y$ can be thought of as the derived category of coherent sheaves 
on a noncommutative K3 surface. Therefore, any smooth moduli space of objects in $\CA_Y$
should be hyperk\"ahler, and the Fano scheme of lines can be realized in this way, see~\cite{KM09}.

The fact that a cubic 4-fold has something in common with a K3 surface can be easily seen from its 
Hodge diamond. In fact, the Hodge diamond of $Y$ is
\begin{equation*}
\begin{smallmatrix}
&&&& 1 \\
&&& 0 && 0 \\
&& 0 && 1 && 0 \\
& 0 && 0 && 0 && 0 \\
0 && 1 && 21 && 1 && 0 \\
& 0 && 0 && 0 && 0 \\
&& 0 && 1 && 0 \\
&&& 0 && 0 \\
&&&& 1 
\end{smallmatrix}
\end{equation*}
and one sees immediately the Hodge diamond of a K3 surface in the primitive part 
of the cohomology of $Y$. There are some other 4-dimensional Fano varieties 
with a similar Hodge diamond. The simplest example is the 4-fold of degree 10 in $\Gr(2,5)$
(an intersection of $\Gr(2,5)$ with a hyperplane and a quadric in $\PP^9$). 
Its Hodge diamond is
\begin{equation*}
\begin{smallmatrix}
&&&& 1 \\
&&& 0 && 0 \\
&& 0 && 1 && 0 \\
& 0 && 0 && 0 && 0 \\
0 && 1 && 22 && 1 && 0 \\
& 0 && 0 && 0 && 0 \\
&& 0 && 1 && 0 \\
&&& 0 && 0 \\
&&&& 1 
\end{smallmatrix}
\end{equation*}
and again its primitive part has K3 type. On a categorical level this follows from the following result

\begin{theorem}[\cite{K09c}]\label{cycat}
Let $X$ be a smooth projective variety of index $m$ with a rectangular Lefschetz decomposition
\begin{equation*}
\BD^b(\coh(X)) = \langle \CB, \CB(1), \dots, \CB(m-1) \rangle
\end{equation*}
of length $m$. Let $Y_d$ be the smooth zero locus of a global section of the line bundle $\CO_X(d)$
for $1 \le d \le m$. Then there is a semiorthogonal decomposition
\begin{equation*}
\BD^b(\coh(Y_d)) = \langle \CA_{Y_d},\CB,\CB(1),\dots,\CB(m-d-1) \rangle
\end{equation*}
and moreover, a power of the Serre functor $\BS_{\CA_{Y_d}}$ is isomorphic to a shift
\begin{equation}
(\BS_{\CA_{Y_d}})^{d/c} = \left[\frac{d\cdot(\dim X + 1) - 2m}c\right],
\qquad
\text{where $c = \gcd(d,m)$}.
\end{equation} 
In particular, if $d$ divides $m$ then
$\CA_{Y_d}$ is a Calabi--Yau category of dimension $\dim X + 1 - 2m/d$.
\end{theorem}

\begin{remark}
Analogously, one can consider a double covering $Y'_d \to X$ ramified in a zero locus
of a global section of the line bundle $\CO_X(2d)$ instead. Then there is an analogous
semiorthogonal decomposition and the Serre functor has the property
\begin{equation}
(\BS_{\CA_{Y'_d}})^{d/c} = \tau^{(m-d)/c}_*\circ\left[\frac{d\cdot(\dim X + 1) - m}c\right],
\end{equation} 
where $\tau$ is the involution of the double covering.
\end{remark}

Applying this result to a 4-fold $Y$ of degree 10 one constructs a semiorthogonal decomposition
$\BD^b(\coh(Y)) = \langle \CA_Y, \CO_Y, \CU_Y^\vee, \CO_Y(1), \CU_Y^\vee(1) \rangle$
with $\CU_Y$ being the restriction of the tautological bundle from the Grassmannian $\Gr(2,5)$
and $\CA_Y$ a Calabi--Yau category of dimension $2$.
Again, for some special $Y$ one can check that $\CA_Y$ is equivalent 
to the derived category of a K3 surface and so altogether we get another family
of noncommutative K3 categories. Moreover, in this case one can also construct
a hyperk\"ahler fourfold from $Y$. One of the ways is to consider the Fano scheme 
of conics on $Y$, see~\cite{IM11}. It turns out that it comes with a morphism to $\PP^5$ with the image
being a singular sextic hypersurface, and the Stein factorization of this map gives
a genus zero fibration over the double covering of the sextic, known as 
a {\sf double EPW sextic}. This is a hyperk\"ahler variety, deformation equivalent 
to the second Hilbert square of a K3 surface.

Finally, there is yet another interesting example. Namely, consider a hyperplane section $Y$
of a 5-fold $X$, which is the zero locus of a global section of the vector bundle
$\Lambda^2\CU_3^\vee \oplus \Lambda^3(W/\CU_3)$ on $\Gr(3,W)$ with $W$ of dimension 7.
This variety $Y$ was found by K\"uchle in \cite{Kuch} (variety $c5$ in his table), and its Hodge 
diamond is as follows
\begin{equation*}
\begin{smallmatrix}
&&&& 1 \\
&&& 0 && 0 \\
&& 0 && 1 && 0 \\
& 0 && 0 && 0 && 0 \\
0 && 1 && 24 && 1 && 0 \\
& 0 && 0 && 0 && 0 \\
&& 0 && 1 && 0 \\
&&& 0 && 0 \\
&&&& 1 
\end{smallmatrix}
\end{equation*}

\begin{conjecture} 
The 5-dimensional variety $X \subset \Gr(3,W)$ has a rectangular Lefschetz decomposition
$\BD^b(\coh(X)) = \langle \CB, \CB(1) \rangle$ with $\CB$
generated by $6$ exceptional objects. 
Consequently, its hyperplane section $Y$ has a semiorthogonal decomposition
$\BD^b(\coh(Y)) = \langle \CA_Y, \CB \rangle$
with $\CA_Y$ being a K$3$ type category.
\end{conjecture}

It would be very interesting to understand the geometry of this variety and to find out,
whether there is a hyperk\"ahler variety associated to it, analogous to the Fano scheme of lines
on a cubic fourfold and the double EPW sextic associated to the 4-fold of degree 10. A natural 
candidate is the moduli space of twisted cubic curves.

\subsection*{Acknowledgements}
I am very grateful to A. Bondal and D. Orlov for their constant help and support.
I would like to thank R.~Abuaf, T.~Bridgeland, T.~Pantev and P.~Sosna
for comments on the preliminary version of this paper.
In my work I was partially supported by RFFI grants NSh-2998.2014.1 and 14-01-00416, 
the grant of the Simons foundation, 
and by AG Laboratory SU-HSE, RF government grant, ag.11.G34.31.0023.

\end{document}